\documentclass[psamsfonts]{article}

\usepackage{amsfonts}
\usepackage{amsmath}
\usepackage{amssymb}
\usepackage{amscd}
\usepackage{amstext}
\usepackage{amsthm}

\title{ Minimal Dynamical Systems and Approximate
Conjugacy
\thanks{Research partially supported by NSF grants DMS 0097903.
         AMS 2000 Subject Classification Numbers:
                         Primary 46L55, 54H20,
                                46L35.
                        Key words:
                        approximate conjugacy, Cantor minimal systems, simple $C^*$-algebras
                                                      \protect\\}}

\author{Huaxin Lin\\
Department of Mathematics\\
University of Oregon\\
Eugene, Oregon 97403-1222\\
and\\
Hiroki Matui\\
Graduate School of Science and Technology\\
Chiba University\\
Chiba 263-8522, Japan}
\date{}

\begin{document}
\maketitle

\renewcommand{\labelenumi}{\textup{(\theenumi)}}
\newcommand{\N}{\mathbb{N}}
\newcommand{\Z}{\mathbb{Z}}
\newcommand{\R}{\mathbb{R}}
\newcommand{\T}{\mathbb{T}}
\newcommand{\C}{\mathbb{C}}
\newcommand{\Aff}{\operatorname{Aff}}
\newcommand{\Inf}{\operatorname{Inf}}
\newcommand{\Homeo}{\operatorname{Homeo}}
\newcommand{\xa}{(X,\alpha)}
\newcommand{\xb}{(X,\beta)}

\newcommand{\CA}{$C^*$-algebra}
\newcommand{\SCA}{$C^*$-subalgebra}
\newcommand{\aue}{approximate unitary equivalence}
\newcommand{\ayue}{approximately unitarily equivalent}
\newcommand{\mops}{mutually orthogonal projections}
\newcommand{\hm}{homomorphism}
\newcommand{\pisca}{purely infinite simple \CA}
\newcommand{\andeqn}{\,\,\,\,\,\, {\rm and} \,\,\,\,\,\,}
\newcommand{\QED}{\rule{1.5mm}{3mm}}
\newcommand{\morp}{contractive completely positive linear map}
\newcommand{\asmorp}{asymptotic morphism}

\newcommand{\arrow}{\rightarrow}
\newcommand{\tdsum}{\widetilde{\oplus}}
\newcommand{\pa}{\|}  
\newcommand{\ep}{\varepsilon}
\newcommand{\id}{{\rm id}}
\newcommand{\aueeps}[1]{\stackrel{#1}{\sim}}
\newcommand{\aeps}[1]{\stackrel{#1}{\approx}}
\newcommand{\dt}{\delta}
\newcommand{\yu}{\fang}
\newcommand{\ca}{{\cal C}_1}
\newcommand{\Ad}{{\rm ad}}

\newtheorem{thm}{Theorem}[section]
\newtheorem{Lem}[thm]{Lemma}
\newtheorem{Prop}[thm]{Proposition}
\newtheorem{Cor}[thm]{Corollary}
\newtheorem{Ex}[thm]{Example}
\newtheorem{Pro}[thm]{Problem}
\newtheorem{NN}[thm]{}
 \theoremstyle{definition}
 \newtheorem{Def}[thm]{Definition}
 \newtheorem{Remark}[thm]{Remark}
\renewcommand{\theequation}{e\,\arabic{section}.\arabic{equation}}

\newcommand{\rforal}{\,\,\,{\rm for\,\,\,all\,\,\,}}
\newcommand{\Ik}{ {\cal I}^{(k)}}
\newcommand{\Iz}{{\cal I}^{(0)}}
\newcommand{\Ii}{{\cal I}^{(1)}}
\newcommand{\Ip}{{\cal I}^{(2)}}
\newcommand{\sless}{{\stackrel{\sim}{<}}}

\begin{abstract}

Several versions of approximate conjugacy for minimal dynamical
systems are introduced. Relation between approximate conjugacy and
corresponding crossed product $C^*$-algebras is discussed. For the
Cantor minimal systems, a complete description is given for these
relations via $K$-theory and $C^*$-algebras. For example, it is
shown that two Cantor minimal systems are approximately
$\tau$-conjugate if and only if they are orbit equivalent and have
the same periodic spectrum. It is also shown that two such systems
are approximately $K$-conjugate if and only if the corresponding
crossed product $C^*$-algebras have the same scaled ordered
$K$-theory. Consequently, two Cantor minimal systems are
approximately $K$-conjugate if and only if the associated
transformation $C^*$-algebras are isomorphic. Incidentally, this
approximate $K$-conjugacy coincides with Giordano, Putnam and
Skau's strong orbit equivalence for the Cantor minimal systems.

\end{abstract}

\section{Introduction}

Let $X$ be a compact metric space and $\alpha, \beta: X\to X$ be
two minimal homeomorphisms. One obtains two dynamical systems $(X,
\alpha)$ and $(X, \beta).$ This paper studies the relation between
these two  systems. The two systems are said to be (flip)
conjugate if there is a homeomorphism $\sigma: X\to X$ such that
$\alpha=\sigma\circ \beta\circ \sigma^{-1}$ (or
$\alpha=\sigma\circ \beta^{-1}\circ \sigma^{-1}$). Let
$C(X)\times_{\alpha}\Z$ be the crossed product \CA\, associated with
the dynamical system $(X, \alpha).$ If $X$ is infinite, the
assumption that $\alpha$ is minimal implies that $C(X)\times
_{\alpha}\Z$ is a unital simple \CA. It follows from a result of
J. Tomiyama (\cite{T2}) (based on an earlier result of M. Boyle)
that $\alpha$ and $\beta$ are flip conjugate if and only if there
exists an isomorphism $\phi: C(X)\times_{\alpha}\Z\to
C(X)\times_{\beta}\Z$ such that $\phi\circ j_{\alpha}=j_{\beta},$
where $j_{\alpha}: C(X)\to C(X)\times_{\alpha}\Z$ and  $ j_{\beta}:
C(X)\to C(X)\times_{\beta}\Z$ are the usual embeddings. While this is
a very interesting result in the interplay between classical
dynamical systems and \CA s, it also suggests that the conjugacy
relation is rather too strong.

Inspired by ergodic theory, Giordano, Putnam and Skau
(\cite{GPS}) introduced the topological version of orbit
equivalences. These notions work nicely when $X$ is the Cantor set.
For example, it is shown in \cite{GPS}, when $X$ is the Cantor
set, then $\alpha$ and $\beta$ are strong orbit equivalent if and
only if the corresponding \CA s $C(X)\times_{\alpha}\Z$ and
$C(X)\times_{\beta}\Z$ are isomorphic. It is a breakthrough  result
which establishes closer relation between \CA\, theory and
dynamical systems in the case that $X$ is the Cantor set. No doubt
it would be more interesting if similar type results could be
established for more traditional topological spaces. It appears
that when $X$ is connected, the situation is rather different. The
recent work of Q. Lin and N. C. Phillips suggests that \CA\,
theory and related $K$-theory may well play an important role in
the study of dynamical systems. We believe that recent development
in the classification theory for simple amenable \CA s via
$K$-theory may also apply to the dynamical systems and $K$-theory
may be proved to be a more useful tool in the study of minimal
dynamical systems.

Motivated by this idea, the first named author has been proposing
a few approximate versions of  conjugacy for more general spaces.
A simple version of approximately conjugate is the following: two
minimal homeomorphisms are approximately conjugate if there exists
a sequence of homeomorphisms $\sigma_n: X\to X$ such that
$f\circ\sigma_n\circ \beta^{-1}\circ \sigma_n^{-1}
\to f\circ\alpha^{-1}$ for
all $f\in C(X).$ When $X$ is not connected, this approximate
version of conjugacy is rather weaker as we will see in Section 4.
When we impose conditions on the conjugate maps $\sigma_n,$ then
the relation becomes more interesting. For example, we introduce
the notion of  approximate $\tau$-conjugacy in which, 
we insist, in addition,  that $\sigma_n$ keep the trace (translation invariant
measure) invariant, and the notion of approximate $K$-conjugacy in
which we require that the conjugate maps $\sigma_n$ induce the
same map on the $K$-theory. We discuss the general relation of
these  approximate conjugacies with corresponding crossed product
$C(X)\times_{\alpha}\Z.$ In the special case that $X$ is the Cantor
set, we show that $\alpha$ and $\beta$ are approximately
$K$-conjugate if and only if $C(X)\times_{\alpha}\Z$ is isomorphic to
$C(X)\times_{\beta}\Z.$ So the approximate $K$-conjugacy is the same
as strong orbit equivalence in this case. We also establish the
relationship between approximate $\tau$-conjugacy and orbit
equivalence in the Cantor set case. These results demonstrate
that perhaps appropriate approximate version of conjugacy is a
right equivalence relation to study and it may bridge the gulf
between the dynamical systems and $C^*$-algebras.

The paper is organized as follows. In Section 2, we list a number
of terminologies used in the paper. We also give a few facts that
are relevant to this paper. In Section 3, we give definitions of
several versions of approximate conjugate relations among minimal
homeomorphisms on compact metric spaces and some consequences. In
Section 4, we discuss weakly approximate conjugate minimal
homeomorphisms on the Cantor set. In Section 5, we show that
approximate conjugate relation is closely related to the orbit
equivalence and strong orbit equivalence in the Cantor minimal systems.
\bigskip

{\bf Acknowledgement} The first named author acknowledges the
support of the National Science Foundation of USA and the Zhi-Jiang
Professorship from East China Normal University. The second named
author is supported by Grant-in-Aid for Young Scientists (B) of
Japan Society for the Promotion of Science. The authors 
benefited from conversations with Takeshi Katsura and N. C.
Phillips.




\section{Preliminaries}

Let $X$ be a compact metric space and $\alpha: X\to X$ be a
minimal homeomorphism. The transformation \CA\, $A_{\alpha}$ may
be identified with the crossed product $C(X)\times_{\alpha}\Z.$
The unitary obtained from $\alpha$ will be denoted by
$u_{\alpha}.$ Denote by $j_{\alpha}: C(X)\to A_{\alpha}$ the
canonical unital embedding.

When $X$ is an infinite set, it is well known that $A_{\alpha}$ is
a simple \CA.

Let $A$ be a unital stably finite \CA. The tracial state space
of $A$ will be denoted by $T(A).$ Let $\Aff(T(A))$ be the
space of all real affine functions on $T(A).$  Denote by $\rho_A: K_0(A)\to
\Aff(T(A))$ the (positive) \hm\, induced by the tracial states.
Let $\tau\in T(A).$ Suppose that $a\in M_k(A)$ for some integer $k,$
we will continue use $\tau$ for the extension $\tau\otimes Tr,$ $Tr$ is
the standard trace (not the normalized one) on $M_k.$

If $h: A\to B$ is a \hm\, from $A$ to another \CA\, $B,$
we denote by $h_{*i}: K_i(A)\to K_i(B)$ ($i=0,1$) the induced
map on $K$-theory. We also use $h_*$ when it is convenient to do so.

If $\sigma: X\to X$ and $f=(f_{ij})_{k\times k}\in M_k(C(X)),$ we
sometimes write  $f\circ\sigma^{-1}$ for
$(f_{ij}\circ\sigma^{-1})_{k\times k}.$ Denote by $\sigma^*:
C(X)\to C(X)$ the induced \hm\, defined by $\sigma^*(f)(x)=
f(\sigma^{-1}(x))$ for $x\in X$ and $f\in C(X).$ 

 We will use
$M_{\infty}(A)$ for the union of $M_k(A),$ $k=1,2,....$

If $A=\lim_{n\to\infty}(A_n,\phi_n)$ is an inductive limit
of \CA s, we use $\phi_{n,\infty}: A_n\to A$ for the induced
\hm\, from $A_n$ to $A.$

\begin{Def}\label{Dasym}
Let $A$ be a separable \CA\, and $B$ another \CA. Let
$\phi_n: A\to B$ be a sequence of \morp s. We say that
$\{\phi_n\}$ is a sequential asymptotic morphism, if
$$
\lim_{n\to\infty}\|\phi_n(ab)-\phi_n(a)\phi_n(b)\|=0\rforal a,
b\in A.
$$
\end{Def}

\begin{Def}\label{DKK}
For each $k,$ we may use $\phi_n: M_k(A)\to M_k(B)$ for
$\phi_n\otimes {\rm id}_{M_k}.$ If $p\in M_k(A)$ is a projection,
then for sufficiently large $n,$ $\phi_n(p)$ is close to a
projection $q\in M_k(B).$ Since close projections are equivalent,
we may write $[\phi_n(p)]=[q].$ In fact, as in 1.8 of \cite{L2} for
example, given any finitely generated subgroup $ G_i\subset K_i(A),$ when $n$
is sufficiently large, $[\phi_n]|_{G_i}$ are well defined
as homomorphisms. We will use this fact whenever it is needed.
In fact, one has the following facts:
Consider $\Phi: A\to \prod_{n\in\N}B$ by defining $\Phi(a)=\{\phi_n(a)\}$
for $a\in A.$ Then $\pi\circ \Phi: A\to \prod_{n\in\N}B/\bigoplus_{n\in\N}B$
is a \hm,
where $\pi: \prod_{n\in \N} B\to \prod_{n\in\N}B/\bigoplus_{n\in\N} B.$
It is well known  (and easy to verify) that $\pi_{*i}:
K_i(\prod_{n\in\N} B)\to K_i(\prod_{n\in\N}B/\bigoplus_{n\in\N} B)$
is surjective and every finitely generated subgroup
of $K_i(\prod_{n\in\N}B/\bigoplus_{n\in\N} B)$ lifts.

\end{Def}

\vspace{0.2in}

\begin{Def}
Let $\xa$ be a Cantor minimal system.
\begin{enumerate}
\item The full group $[\alpha]$ of $\xa$ is the subgroup of
all homeomorphisms from $X$ to $X$  that preserves
every orbit of $\alpha$.
To any $\gamma\in[\alpha]$ is associated a map $n:X\rightarrow \Z$,
defined by $\gamma(x)=\alpha^{n(x)}(x)$ for $x\in X$.
\item The topological full group $[[\alpha]]$ of $\xa$ is
the subgroup of all homeomorphisms $\gamma\in[\alpha]$,
whose associated map $n:X\rightarrow \Z$ is continuous.
\end{enumerate}
\end{Def}

%

\begin{Lem}\label{Llift1}
Let $X$ be the Cantor set and $\alpha: X\to X$ be a minimal
homeomorphism.
Let $p\in C(X)$ be a projection, $x\in K_0(A_{\alpha})_+\setminus\{0\}$
such that $x\le (j_{\alpha})_{*0}(p).$ then there
is a projection $q\le p$ such that $(j_{\alpha})_{*0}(q)=x.$
\end{Lem}
\begin{proof}
Let $Y=\{y\}$ be a subset of $X$ consisting of a single point. It
follows from Corollary 1 in \cite{GPS} that $i_*:
K_0(A_Y^{\alpha})\to K_0(A_\alpha)$ is an isomorphism of ordered
groups preserving the order unit. Thus there is a non-zero
projection $q\in A_Y^{\alpha}$ such that $i_*([q])=x.$ It follows
from Theorem VIII.6.2 and its proof in \cite{Dk} that
$(j_{\alpha})_*$ maps $K_0(C(X))_1^+$ onto $K_0(A_Y^\alpha)_1^+$
and that there exists a projection $q'\le p$ such that
$[(j_{\alpha})_{*0}(q')]=[q].$
\end{proof}

\begin{Lem}\label{Lext3}
Let $X$ be the Cantor set and $\alpha: X\to X$ be a minimal
homeomorphism.
Suppose that $x_1,x_2,...,x_m\in K_0(A)_+\setminus\{0\}$
satisfy $\sum_{i=1}^mx_i=[1_A].$
Then there are mutually orthogonal projections $q_1,q_2,...,q_m\in
C(X)$ such that $\sum_{i=1}^mq_i=1_{C(X)}$ and
$$
(j_{\alpha})_{*0}([q_i])=x_i,\,i=1,2,...,m.
$$
\end{Lem}
\begin{proof}
In the same way as Lemma \ref{Llift1}, there exists
a projection $q_1$ such that $(j_{\alpha})_{*0}([q_1])=x_1$.
By applying Lemma \ref{Llift1} to $x_2$ and $1-q_1$,
we get $q_2$.
The projections $q_3,q_4,\dots,q_{m-1}$ are obtained in the same way.
Put $q_m=1-(q_1+q_2+\dots +q_{m-1})$.
\end{proof}
\bigskip

\begin{thm}\label{PreT}
Let $X$ be the Cantor set and $\alpha: X\to X$ be a minimal
homeomorphism. Suppose that $\phi: C(X,\Z)\to K_0(A_{\alpha})$ is
a  positive \hm\, for which $\phi(1_{C(X)})=[1_{A_{\alpha}}].$
\begin{enumerate}
\item Then there exists a unital positive \hm\, $\psi: C(X,\Z)\to C(X,\Z)$ such
that
$$
(j_{\alpha})_{*0}\circ \psi=\phi.
$$
\item If $\phi(C(X,\Z)_+\setminus \{0\})\subset K_0(A_{\alpha})_+\setminus
\{0\},$ then one can require that $\psi$ to be injective.
\item Suppose that $\beta: X\to X$ is another minimal homeomorphism,
$\gamma: K_0(A_{\beta})\to K_0(A_{\alpha})$ is a unital
order isomorphism and $\phi=\gamma\circ (j_{\beta})_{*0},$
then one can require
that $\psi$ has also an invertible  positive \hm.
\end{enumerate}

\end{thm}

\begin{proof}
Write $C(X, \Z)=\lim_{n\to\infty} (\Z^{k(n)}, \zeta_n).$ Here we
use the usual order on $\Z^{k(n)}.$  Each $\zeta_n$ is injective, unital  and
positive.
We regard $\Z^{k(n)}$ as the subgroup of $G$.
Note that $C(X,\Z)_+=\bigcup_{n=1}^{\infty}\Z^{k(n)}_+.$

Let $e^{(n)}_1=(1,0,...,0), e^{(n)}_2=(0,1,0,...,0),...,
e^{(n)}_{k(n)}=(0,...,0,1),$ $n=1,2,....$

Fix $n,$  let $x(n,1), x(n,2),...,x(n,k(n))\in K_0(A_{\alpha})_+$ such
that $\phi(e^{(n)}_i)=x(n,i),$ $i=1,2,..., k(n)$ (some of them
could be zero).
By applying Lemma \ref{Lext3} above, we obtain a positive \hm\,
$h_n: \Z^{k(n)}\to C(X,\Z)$ such that
$$
(j_{\alpha})_{*0}\circ h_n=
\phi|_{\Z^{(k(n)}}.
$$

Moreover, if $\phi(C(X,\Z)_+\setminus \{0\})\subset
K_0(A_{\alpha})_+\setminus \{0\},$ $x(n,i)$ above are non-zero and
by Lemma \ref{Lext3}, we may assume that
$h_n(e^{(n)}_i)$ are mutually orthogonal projections
in $C(X)$ (Here we regard $C(X,\Z)$ as a subset of $C(X)$).
In particular, in this case, $h_n$ is injective.

We note that in
$\Z^{k(n)},$ $1_{C(X)}=\sum_{i=1}^{k(n)}e_i^{(n)}.$ Furthermore,
we may write
$$
\sum_{j\in S(n,i)}e_j^{(n+1)}=e_i^{(n)},
$$
where $\{S(n,1), S(n,2),...,S(n,k(n))\}$ are disjoint and
$\bigcup_iS(n,i)=\{1,2,...,k(n+1)\}.$ Let
$p_i=h_n(e_i^{(n)}).$
By applying Lemma \ref{Llift1} and using the same
argument used in Lemma \ref{Lext3},
we obtain projections $p(n+1,j)\in
C(X,\Z)_+$  such that
$$
\sum_{j\in S(n,i)}p(n+1,j)=p_i=h_n(e_i^{(n)}) \andeqn
(j_{\alpha})_{*0}(p(n+1,j))=\phi(e_j^{(n+1)}).
$$
Define $h_{n+1}: \Z^{k(n+1)}\to C(X,\Z)$ by
$h_{n+1}(e_j^{(n+1)})=p(n+1,j),$ $j=1,2,...,k(n+1).$ Clearly
$h_{n+1}(1_{C(X)})=1_{C(X)}.$ More importantly,
$$
h_{n+1}|_{\Z^{k(n)}}=h_n,n=1,2,....
$$
Thus $\{h_n\}$ defines a positive \hm\, $\psi: C(X,\Z)\to C(X,\Z)$
which preserves $1_{C(X)}.$ Moreover,
$$
(j_{\alpha})_{*0}\circ \psi=
\phi.
$$

Note that in the case (2), each $h_n$ is
injective, it follows that $\psi$ is also injective.
This proves (1) and (2).

To see (3) holds, we
note that, since both $K_0(A_{\alpha})$ and $K_0(A_{\beta})$ are
simple dimension groups,\\
$\phi(C(X,\Z)_+\setminus \{0\})\subset
K_0(A_{\alpha})_+\setminus \{0\}.$

Suppose that $h_n$ has been constructed.
We may assume that $h_n(e^{(n)}_i)\subset \Z^{k(n+1)},$ without
loss of generality.

We write $\sum_{j\in F(n,i)}e^{(n+1)}_j=h_n(e_i^{(n)}),$
where $\{F(n,i)\}$ are mutually disjoint and \\
$\bigcup_jF(n,j)=\{1,2,...,k(n+1)\}.$
Suppose that $y(n,j)=\gamma^{-1}\circ (j_{\alpha})_{*0}(e^{(n+1)}_j),$
$j\in F(n,i).$
Note
that $\gamma\circ (j_{\beta})_{*0}(e_i^{(n)})=(j_{\alpha})_{*0}\circ
h_n(e_i^{(n)}),$ $i=1,2,...,k(n).$
Therefore
$$
(j_{\beta})_{*0}(e_i^{(n)})=\sum_{j\in F(n,i)}y(n,j).
$$
By Lemma \ref{Llift1} and the proof of Lemma \ref{Lext3},
we obtain mutually orthogonal non-zero projections\\
$p_j$ ($j\in F(n,i)$) such that
$$
\sum_{j\in F(n,i)}p_j=e_i^{(n)}\andeqn
(j_{\beta})_{*0}(p_i)=
y(n,j)
$$
for some $m>n.$
Without loss of generality, we may assume
that $m=n+2.$

By defining $H_n(e^{(n+1)}_j)=p_j,$ we obtain a unital positive injective
\hm\, $H_n: \Z^{k(n+1)}\to \Z^{k(n+2)}\subset C(X,\Z)$ such that
$$
(j_{\beta})_{*0}\circ H_n=\gamma^{-1}\circ (j_{\alpha})_{*0}|_{\Z^{k(n+1)}}.
$$
Moreover, $H_n\circ h_n$ is the identity on $\Z^{k(n)}$.
In other words, we have the following
commutative diagram.
$$
\begin{array}{ccc}
\Z^{k(n)} & \stackrel{\zeta_{n, n+2}}{\longrightarrow} & \Z^{k(n+2)}\\
\downarrow_{h_n} & \nearrow_{H_n} \\
\Z^{k(n+1)}&
\end{array}
$$

Exactly the same argument allows us to construct $h_{n+1}:
\Z^{k(n+2)}\to Z^{k(n+3)}$ (by perhaps choosing a number larger
than $n+3$) such that $(j_{\alpha})_{*0}\circ h_{n+1}=\gamma\circ
j_{\beta}|_{\Z^{k(n+2)}}$ and
we have the following commutative diagram.
$$
\begin{array}{ccc}
\Z^{k(n)} & \stackrel{\zeta_{n, n+2}}{\longrightarrow} & \Z^{k(n+2)}\\
\downarrow_{h_n} & \nearrow_{H_n} & \downarrow_{h_{n+1}}\\
\Z^{k(n+1)}& \stackrel{\zeta_{n+1,n+3}}{\longrightarrow} & Z^{k(n+3)}\\
\end{array}
$$

Continuing in this fashion, we obtain $h_{n+m}: \Z^{k(n+2m)}\to
\Z^{k(n+2m+1)}$ and $H_{n+m}: \Z^{k(n+2m+1)} \to \Z^{k(n+2m+2)}$ such
that $(j_{\alpha})_{*0}\circ h_{n+m}=\gamma\circ
(j_{\beta})_{*0}|_{\Z^{k(n+2m)}}$ and $(j_{\beta})_{*0}\circ
H_{n+m}=\gamma^{-1}\circ (j_{\alpha})_{*0}|_{\Z^{(n+2m+1)}}.$
Furthermore, we have the following intertwining:
$$
\begin{array}{cccccc}
\Z^{k(n)} & \stackrel{\zeta_{n, n+2}}{\longrightarrow} & \Z^{k(n+2)}
&  \stackrel{\zeta_{n+2,n+4}}{\longrightarrow} \cdots &
\to\cdots & C(X, \Z)\\
\downarrow_{h_n} & \nearrow_{H_n} & \downarrow_{h_{n+1}}
& \nearrow_{H_{n+1}} &\\
\Z^{k(n+1)}& \stackrel{\zeta_{n+1,n+3}}{\longrightarrow} & \Z^{k(n+3)}
&  \stackrel{\zeta_{n+3, n+5}}{\longrightarrow} \cdots &
\to \cdots &C(X,\Z) \\
\end{array}
$$

This implies that
we have  two invertible unital order \hm s $\psi, \psi_1: C(X, \Z)
\to C(X,\Z)$  such that $\psi_1=\psi^{-1}.$
By the restriction on $h_n$ and $H_n,$ we also have
$$
(j_{\alpha})_{*0}\circ \psi=\gamma\andeqn (j_{\beta})_{*0}\circ \psi^{-1}.
$$
\end{proof}


\begin{Lem}\label{ILeffchoi}
Let $A$ be a separable simple unital \CA\, and $B$ be another \CA.
Suppose that $\{\phi_n\}: A\to B$ be a (unital) sequential
asymptotic morphism. Then, for any $a\not=0$ in $A,$
$$
\lim_{n\to\infty}\|\phi_n(a)\|=\|a\|.
$$
\end{Lem}

\begin{proof}
Note that $\|\phi_n(a)\|\le \|a\|.$ Suppose that there is a
nonzero element $a\in A$ such that
$$
\liminf_n\|\phi_n(a)\|<\|a\|.
$$
Then there is  a subsequence $\{n(k)\}$ such that
\begin{eqnarray}\label{eeffchoi}
\lim_{k\to\infty}\|\phi_{n(k)}(a)\|<\|a\|.
\end{eqnarray}
Let $\Phi: A\to \prod_{k=1}^{\infty} B$ be defined by
$\Phi(a)=\{\phi_{n(k)}(a)\}_k.$ Let $\Pi: \prod_{k=1}^{\infty}B\to
\prod_{k=1}^{\infty}B/\bigoplus_{k=1}^{\infty}B$ be the quotient
map. Then $\Pi\circ \Phi:
A\to\prod_{k=1}^{\infty}B/\bigoplus_{k=1}^{\infty}B$ is a \hm.
Since $A$ is simple, it is a monomorphism. In particular,
$\|\Pi\circ\Phi(a)\|=\|a\|.$ This contradicts (\ref{eeffchoi}).
\end{proof}

\section{ Approximate conjugacy}

Let $X$ be a compact metric space and $\alpha, \beta: X\to X$ be
two minimal homeomorphisms.
Recall that $\alpha$ and $\beta$ are (flip) conjugate if
there is a homeomorphism $\sigma: X\to X$ such that
$\sigma\circ \alpha\circ \sigma^{-1}=\beta$
(or $\sigma\circ \alpha\circ \sigma^{-1}=\beta^{-1}).$
It is in general rather strong equivalence relation.
We propose a few approximate version of conjugacy.

We start with the following:

\begin{Def}\label{IIID1}
Let $X$ be a compact metric space and $\alpha, \beta: X\to X$
be two homeomorphisms. We say $\alpha$ and $\beta$ are {\it weakly
approximately conjugate} if there are homeomorphisms $\sigma_n,
\gamma_n: X\to X$ such that
\[
\lim_{n\to\infty}(\sigma_n\circ \alpha\circ \sigma_n^{-1})^*(f)=\beta^*(f)
\andeqn \lim_{n\to\infty} (\gamma_n\circ \beta\circ \gamma_n^{-1})^*(f)=
\alpha^*(f)
\]
for all $f\in C(X).$
\end{Def}

\begin{Prop}\label{wac>asym}
Let $X$ be a compact metrizable space and let $\alpha$ and $\beta$
be minimal homeomorphisms on $X$. Let
$A=C(X)\times_{\alpha}{\mathbb Z}$ and
$B=C(X)\times_{\beta}{\mathbb Z}.$  If there exist
homeomorphisms $\sigma_n:X\to X$ such that
$(\sigma_n\circ \alpha\circ \sigma_n^{-1})^*(f)
\rightarrow (\beta)^*(f)$ for
every $f\in C(X),$  then there exists a unital  asymptotic
morphism $\{\psi_n\}$ from $B$ to $A$ such that
\[ \lim_{n\rightarrow\infty}
||\psi_n(u_\beta)-u_\alpha||=0 \]
and
\[ \lim_{n\rightarrow\infty}
||\psi_n(j_\beta(f))-j_\alpha(f\circ\sigma_n)||=0 \] for all
$f\in C(X)$. In particular, for every $b\in B$,
\[ \lim_{n\rightarrow\infty}
||\psi_n(E_B(b))-E_A(\psi_n(b))||=0,\]
where $E_A$ and $E_B$ are the canonical conditional expectations
onto $C(X)$ defined on $A$ and $B$.
\end{Prop}
\begin{proof}
Let us denote the quotient map from $\prod_{n\in\N}A$ to
$\prod_{n\in\N}A/\bigoplus_{n\in\N}A$ by $\pi$.
Define $\Psi(j_\beta(f))=\pi((j_\alpha(f\circ\sigma_n))_n)$ for
$f\in C(X)$
and $\Psi(u_\beta)=\pi((u_\alpha)_n).$  One checks that $\Phi$
gives   a homomorphism
from $B$ to $\prod_{n\in\N}A/\bigoplus_{n\in\N}A$,
because
\[ u_\alpha j_\alpha(f\circ\sigma_n)u_\alpha^*
-j_\alpha(f\circ\beta^{-1}\sigma_n)
=j_\alpha(f\circ\sigma_n\alpha^{-1}-f\circ\beta^{-1}\sigma_n)
\rightarrow 0. \]
Since $B$ is
nuclear, there exists a unital completely positive lifting of
$\Psi$ (\cite{CE}), that is, we can find a sequence $\{\psi_n\}$
of unital completely positive linear maps from $B$ to $A$ such that
$$
\lim_{n\to\infty}
\left\|\psi_n\left(\sum_{k=-N}^Lf_ku_{\beta}^k\right)-
\sum_{k=-N}^L(f_k\circ\sigma_n)u_{\alpha}^k\right\|=0
$$
for all $\sum_{k=-N}^Lf_ku_{\beta}^k\in B$.
Because $\pi\circ \Phi$ is a homomorphism, $\{\psi_n\}$ is
an asymptotic morphism.
\end{proof}

Let $Q_A=\prod_{n\in \N} A/\bigoplus_{n\in\N} A,$
$Q_{C(X)}=\prod_{n\in \N}C(X)/\bigoplus_{n\in\N} C(X).$
Denote by ${\bar \alpha}$ the induced action on $Q_{C(X)}.$
Let $Q_{\alpha}=Q_{C(X)}\times_{{\bar\alpha}}\Z.$
It is easy to see that there is an obvious unital embedding from
$Q_{\alpha}$ to $Q_A.$

\begin{Cor}\label{IIIC1}
In the situation of Proposition \ref{wac>asym}, one has
the following commutative diagram:
\[
{\small
\put(-190,0){$K_0(C(X))$}
\put(-30,0){$K_0(C(X))$} \put(130,0){$K_0(B)$}
\put(-115,-40){$K_0(Q_{C(X)})$} \put(-30,-40){$K_0(Q_{C(X)})$}
\put(55,-40){$K_0(Q_\alpha)$} \put(-115, -70){$K_1(Q_\alpha)$}
\put(-30,-70){$K_1(Q_{C(X)})$} \put(55,-70){$K_1(Q_{C(X)} )$}
\put(-190,-110){$K_1(B)$} \put(-30,-110){$K_1(C(X))$}
\put(130,-110){$K_1(C(X))$}
\put(-140, 2){\vector(1,0){105}}
\put(20,1){\vector(1,0){105}} \put(-153,-7){\vector(1,-1){30}}
\put(-10,-3){\vector(0,-1){25}} \put(125,-7){\vector(-1,-1){25}}
\put(-69,-38){\vector(1,0){35}} \put(19,-38){\vector(1,0){32}}
\put(-172, -97){\vector(0,1){90}} \put(-105,-60){\vector(0,1){15}}
\put(70, -45){\vector(0,-1){15}} \put(140,-7){\vector(0,-1){90}}
\put(-37,-68){\vector(-1,0){35}} \put(50,-68){\vector(-1,0){32}}
\put(-153,-102){\vector(1,1){30}}
\put(125,-100){\vector(-1,1){25}} \put(-10,-100){\vector(0,1){25}}
\put(-35,-108){\vector(-1,0){105}} \put(125,-108){\vector(-1,0){105}}
\put(-142,-15){$(\pi\circ{\tilde\sigma})_*$}
\put(-42,-15){$(\pi\circ{\tilde\sigma})_*$}
\put(105, -15){$\Psi_*$}
\put(-155,-92){$\Psi_*$} \put(-42, -92){$(\pi\circ{\tilde\sigma})_*$}
\put(80, -92){$(\pi\circ{\tilde \sigma})_*$}
\put(-105,4){${\rm id}_*-\beta_*$} \put(67,4){$(j_{\beta})_*$}
\put(45,-104){${\rm id}_*-\beta_*$} \put(-100,-104){$(j_{\beta})_*$}
},\]
where ${\tilde\sigma}: C(X)\to \prod_{n\in\N}C(X)$ is defined by
${\tilde\sigma}(f)=(f\circ \sigma_n)_n,$
the map from $K_i(Q_{C(X)})$ to $K_i(Q_{C(X)})$ is
$\id-{\bar\alpha}_*$ and the map from $K_i(Q_{C(X)})$ to
$K_i(Q_\alpha)$ is induced by the natural inclusion.

Moreover, we also have the following commutative diagram:
\[
{\small
\put(-190,0){$K_0(C(X))$}
\put(-30,0){$K_0(C(X))$} \put(130,0){$K_0(B)$}
\put(-115,-40){$Q_{K_0(C(X))}$}
\put(-30,-40){$Q_{K_0(C(X))}$}
\put(55,-40){$Q_{K_0(A)}$} \put(-115, -70){$Q_{K_1(A)}$}
\put(-30,-70){$Q_{K_1(C(X))}$} \put(55,-70){$Q_{K_1(C(X))}$}
\put(-190,-110){$K_1(B)$} \put(-30,-110){$K_1(C(X))$}
\put(130,-110){$K_1(C(X))$}
\put(-140, 2){\vector(1,0){105}}
\put(20,1){\vector(1,0){105}} \put(-153,-7){\vector(1,-1){30}}
\put(-10,-3){\vector(0,-1){25}} \put(125,-7){\vector(-1,-1){25}}
\put(-69,-38){\vector(1,0){35}} \put(19,-38){\vector(1,0){32}}
\put(-172, -97){\vector(0,1){90}} \put(-105,-60){\vector(0,1){15}}
\put(70, -45){\vector(0,-1){15}} \put(140,-7){\vector(0,-1){90}}
\put(-37,-68){\vector(-1,0){35}} \put(50,-68){\vector(-1,0){32}}
\put(-153,-102){\vector(1,1){30}}
\put(125,-100){\vector(-1,1){25}} \put(-10,-100){\vector(0,1){25}}
\put(-35,-108){\vector(-1,0){105}} \put(125,-108){\vector(-1,0){105}}
\put(-142,-15){$\theta_0\circ(\pi\circ{\tilde\sigma})_*$}
\put(-2,-15){$\theta_0\circ(\pi\circ{\tilde\sigma})_*$}
\put(-60, -92){$\theta_1\circ(\pi\circ{\tilde\sigma})_*$}
\put(62, -92){$\theta_1\circ(\pi\circ{\tilde \sigma})_*$}
\put(-105,4){${\rm id}_*-\beta_*$} \put(67,4){$(j_{\beta})_*$}
\put(45,-104){${\rm id}_*-\beta_*$} \put(-100,-104){$(j_{\beta})_*$}
},\]
where $Q_{K_i(C(X))}$ denotes
$\prod_{n\in\N}K_i(C(X))/\bigoplus_{n\in\N} K_i(C(X))$,
$Q_{K_i(A)}$ denotes
$\prod_{n\in\N}K_i(A)/\bigoplus_{n\in\N} K_i(A)$
and $\theta_*$ is the natural homomorphisms
from $K_*(Q_{C(X)})$ to $Q_{K_*(C(X))}$.
\end{Cor}

\begin{proof}
By Proposition \ref{wac>asym}, $\pi\circ\tilde{\sigma}$ induces
a homomorphism from $B$ to $Q_{\alpha}$.
Since the Pimsner-Voiculescu exact sequences behave naturally
for a covariant homomorphism,
the commutativity of the first diagram follows.

We will show the commutativity of the second diagram.
Let us consider the following diagram:
\[ \begin{CD}
K_*(Q_{C(X)}) @>\id-\alpha_*>> K_*(Q_{C(X)}) \\
@V\theta_*VV @VV\theta_*V \\
Q_{K_*(C(X))} @>>> Q_{K_*(C(X))},
\end{CD} \]
where the horizontal map below is defined by
\[ (x_n)_n+\bigoplus_{n\in\N}K_*(C(X))\mapsto
(x_n-\alpha_*(x_n))_n+\bigoplus_{n\in\N}K_*(C(X)). \]
Then the commutativity is clear.

Next, let us consider
\[ \begin{CD}
K_*(Q_{C(X)}) @>>> K_*(Q_\alpha) \\
@V\theta_*VV @VVV \\
Q_{K_*(C(X))} @>>> Q_{K_*(A)},
\end{CD} \]
where the top horizontal map is induced by the natural inclusion,
the lower horizontal map  is induced by the natural map
$K_*(C(X))\rightarrow K_*(A)$ and
the vertical map of the right-hand side is the composition
of the natural map from $K_*(Q_\alpha)$ to $K_*(Q_A)$ and
the natural map from $K_*(Q_A)$ to $Q_{K_*(A)}$.
Then, again, the commutativity is easily verified.

Finally, we consider the diagram:
\[ \begin{CD}
K_*(Q_\alpha) @>>> K_{1-*}(Q_{C(X)}) \\
@VVV @VVV \\
Q_{K_*(A)} @>>> Q_{K_{1-*}(C(X))}.
\end{CD} \]
The vertical maps and the top horizontal map
have been already defined. The lower horizontal map is induced by
the map from $K_*(A)$ to $K_{1-*}(C(X))$
appeared in Pimsner-Voiculescu exact sequence.
Therefore we can check the commutativity.
\end{proof}


%

\begin{Def}\label{IIIDpro}
Let $X$ be a compact metric space and $\alpha: X\to X$ be a
minimal homeomorphism. A positive element $f\in M_k(C(X))$ is
said to be $\alpha$-projection valued, if there is a projection
$p=\sum_{i\in {\Z}}f_nU_{\alpha}^n\in M_k(A_{\alpha})$ such that
$f_0=f,$ where $U_{\alpha}={\rm
diag}(u_{\alpha},u_{\alpha},...,u_{\alpha}).$ 

\end{Def}

\begin{Prop}\label{IIIPpv}
If $f\in M_n(C(X))$ is $\alpha$-projection valued, then there is a
projection $p\in M_n(A_{\alpha})$ such that $\tau(p)=\tau(f)$ for
every $\alpha$-invariant tracial state $\tau$.

If furthermore,
$(j_{\alpha})_{*0}(K_0(C(X))_+)\supset K_0(A_{\alpha})_+,$ then
there is a projection $q\in M_k(C(X))$ (for some integer $k$) such
that $\tau(f)=\tau(q)$ for all $\alpha$-invariant tracial states
$\tau$.
\end{Prop}

\begin{proof}
Suppose that $p=\sum_{j\in \Z}f_jU_{\alpha}^j$ for which $f_0=f.$
Let $\tau$ be the tracial state on $A_{\alpha}$ and $E$ be the
canonical faithful expectation from $A_{\alpha}$ onto $j_\alpha(C(X))$.
Then $E(p)=f_0=f$
and, consequently, $\tau(p)=\tau(f).$

To see the second part of the statement, we let $p=\sum_{j\in \Z}
f_jU_{\alpha}^j$ be a projection in $M_n(A_{\alpha})$ with
$f_0=f.$ Since  $(j_{\alpha})_{*0}: K_0(C(X))_+\to
K_0(A_{\alpha})_+$ is surjective, there is a projection $q\in
M_k(C(X))$ for some integer $k>0$ such that
$(j_{\alpha})_{*0}([q])=[p].$ It follows that there is an integer
$m>0$ and a partial isometry $v\in M_{k+m}(A_{\alpha})$ such that
$$
v^*v=p\oplus 1_m\andeqn vv^*=q\oplus 1_m.
$$
It follows that $\tau(p)+m=\tau(q)+m$ for all $\tau\in T(A_{\alpha}).$
Thus
$\tau(p)=\tau(q)$ for all $\tau\in T(A_{\alpha}).$
\end{proof}


Comparing with the definition of conjugacy, one may impose some
restriction on the maps $\sigma_n$ so that $\sigma_n$ and
$\sigma_{n+1}$ may have something in common. In the following
definition, we require that the conjugate maps $\sigma_n$ keep the
invariant measures unchanged. Later  in Definition \ref{DKC}, we will give
another version of approximate conjugacy where we require that the
conjugate maps preserve the $K$-theory.

\begin{Def}\label{IIIDtauc}

We say that  $\alpha$ and $\beta$ are {\it approximately
$\tau$-conjugate}, if $\alpha$ and $\beta$ are weakly
approximately conjugate, and, if in addition, for any  projection
$f\in M_k(C(X))$ there exists $N$ such that
$$
\tau'(f\circ\gamma_n)=\tau'(f\circ\gamma_{n+1}), \ \
\tau(f\circ\sigma_n)=\tau(f\circ\sigma_{n+1}),
$$
$$
\tau(f\circ\gamma_n\circ\sigma_m)=\tau(f)
\andeqn \tau'(f\circ\sigma_n\circ\gamma_m)=\tau'(f)
$$
for every $\alpha$-invariant measure $\tau,$ $\beta$-invariant
measure $\tau'$ and all $m\ge n\ge N.$

Moreover,
for any $\alpha$-projection valued functions $f\in M_\infty(C(X))$
and $\beta$-projection valued function $g\in M_\infty(C(X)),$
there exist $\beta$-projection valued function $f'\in M_\infty(C(X))$
and $\alpha$-projection valued function $g'\in M_\infty(C(X))$
such that
$$
\lim_{n\to\infty}\tau'(f\circ\gamma_n)=\tau'(f'), \ \
\lim_{n\to\infty}\tau(g\circ\sigma_n)=\tau(g')
\andeqn
$$
$$
\lim_{n\to\infty}\limsup_{m\to\infty}|
\tau(f\circ\gamma_n\circ\sigma_m)-\tau(f)|=0,
\lim_{n\to\infty}\limsup_{m\to\infty}|
\tau'(g\circ\sigma_n\circ\gamma_m)-\tau'(g)|=0
$$
uniformly for
$\alpha$-invariant normalized measures $\tau$
and $\beta$-invariant normalized measures $\tau'$.
\end{Def}

\begin{Prop}\label{IIIP2}
Let $X$ be a compact metric space, $\alpha, \beta: X\to X$ be two
minimal homeomorphisms. Let $A=A_{\alpha}$ and $B=A_{\beta}.$
Suppose that $\alpha$ and $\beta$ are approximately
$\tau$-conjugate. Then there exist two unital sequential
asymptotic morphisms $\{\phi_n\}: A \to B$ and $\{\psi_n\}: B\to
A$ such that, for each projection $p\in M_{\infty}(A)$ and
projection $q\in M_{\infty}(B),$ there exist projections $p'\in
M_{\infty}(B)$ and $q'\in M_{\infty}(A)$ such that
\begin{eqnarray}\label{etreq}
\lim_{n\to\infty}\tau'(\phi_n(p))=\tau'(p'), \ \
\lim_{n\to\infty}\tau(\psi_n(q))=\tau(q')\andeqn
\end{eqnarray}
\begin{eqnarray}\label{etr=1}
\lim_{n\to\infty}\limsup_{m\rightarrow\infty}
|\tau(\psi_m\circ \phi_n(p))-\tau(p)|=0, \ \
\lim_{n\to\infty}\limsup_{m\rightarrow\infty}
|\tau'(\phi_m\circ \psi_n(q))-\tau'(q)|=0
\end{eqnarray}
uniformly for $\tau\in T(A)$ and $\tau'\in T(B)$.
\end{Prop}

\begin{proof}
Proposition \ref{wac>asym} applies and yields asymptotic morphisms
$\{\phi_n\}$ and $\{\psi_n\}$.
Let $E_A$ (resp. $E_B$) be the canonical conditional expectation
from $A$ (resp. $B$) to $C(X)$. Let $p\in M_\infty(A)$ be a
projection. Then $E_A(p)$ is an $\alpha$-projection valued
function. By Definition \ref{IIIDtauc} and \ref{IIIDpro}, there is
a projection $q\in M_\infty(B)$ such that
\[ \lim_{n\rightarrow\infty}
\tau'(E_A(p)\circ\gamma_n)=\tau'(E_B(q))=\tau'(q) \]
uniformly for $\tau'\in T(B)$. Since
\begin{align*}
&||E_A(p)\circ\gamma_n-E_B(\phi_n(p))|| \\
&=||E_A(p)\circ\gamma_n-\phi_n(E_A(p))||+
||\phi_n(E_A(p))-E_B(\phi_n(p))||\rightarrow0,
\end{align*}
we get the desired equality.

Moreover, by Definition \ref{IIIDtauc} and Proposition \ref{wac>asym},
for any $\ep>0$,
there is $N\in\N$ such that
\[ \sup_{\tau\in T(A)}\limsup_{m\to\infty}
|\tau(E_A(p)\circ\gamma_n\circ\sigma_m)-\tau(p)|<\ep \]
and
\[ ||E_A(p)\circ\gamma_n-E_B(\phi_n(p))||<\ep,\]
for all $n\geq N$. Then we have
\begin{align*}
& \max_{\tau\in T(A)}\limsup_{m\rightarrow\infty}
|\tau(\psi_m(\phi_n(p)))-\tau(p)| \\
&=\max_{\tau\in T(A)}\limsup_{m\rightarrow\infty}
|\tau(E_A(\psi_m(\phi_n(p))))-\tau(p)| \\
&=\max_{\tau\in T(A)}\limsup_{m\rightarrow\infty}
|\tau(\psi_m(E_B(\phi_n(p))))-\tau(p)| \\
&\leq\max_{\tau\in T(A)}\limsup_{m\rightarrow\infty}
|\tau(\psi_m(E_A(p)\circ\gamma_n))-\tau(p)|+\ep \\
&=\max_{\tau\in T(A)}\limsup_{m\rightarrow\infty}
|\tau(E_A(p)\circ\gamma_n\circ\sigma_m)-\tau(p)|+\ep<2\ep
\end{align*}
for all $n\geq N$.
The other equalities are shown in the same way.
\end{proof}

\begin{thm}\label{IIITtr}
Let $X$ be a compact metric space, $\alpha,\beta: X\to X$ be
minimal homeomorphisms. Suppose that $\alpha$ and $\beta$ are
approximately  $\tau$-conjugate. Then there is an order unit
preserving order isomorphism between $\rho_A(K_0(A_{\alpha}))$ and
$\rho_B(K_0(A_{\beta})).$
\end{thm}

\begin{proof}
Let $A=A_{\alpha}$ and $B=A_{\beta}.$

Let $\{\phi_n\}: A\to B$ and $\{\psi_n\}$ be unital sequential asymptotic
morphisms as described in Proposition \ref{IIIP2}.

To see that $\rho_A(K_0(B))$ is order isomorphic to $\rho_B(K_0(B)),$
we define a map $H:\rho_A(K_0(A))\to
\rho_B(K_0(B))$ by
$$
H(\hat{p})(\tau')=\lim_{n\to\infty} \tau'\circ \phi_n(p)
$$
for each projection $p\in M_\infty(A)$ and $\tau'\in T(B).$

To show that $H$ is a well defined map from $\rho_A(K_0(A))$ to
$\Aff(T(B)),$ we first need to verify that if $q$ is another
projection such that $\tau(p)=\tau(q)$ for all $\tau\in T(A),$
then $H(\hat{p})=H(\hat{q}).$

Let $\tau'\in T(B).$ Define $s_n=\tau'\circ\phi_n.$ Then $s_n$ is a
state of the unital \CA\, $A.$ Let $s$ be a weak limit of
$\{s_n\}.$ Then it is easy to check that $s$ gives a tracial state
of $A.$ It follows that $s(p)=s(q).$ Thus
$$
\lim_{n\to\infty} (s_n(p)-s_n(q))=0.
$$
This shows that $H(\hat{p})=H(\hat{q}).$
By (\ref{etreq}), we see that the image of $H$ is in $\rho_B(K_0(B)).$
In other words, $H$ maps
$\rho_A(K_0(A))$ to $\rho_B(K_0(B)).$ It is then easy to check
that $H$ is a positive \hm\, and preserves the order unit.
Similarly,
$$
H'(\hat{q})(\tau)=\lim_{n\to\infty} \tau\circ\psi_n(q)
$$
defines a positive \hm\, from $\rho_B(K_0(B))$ to $\rho_A(K_0(A))$
which preserves the order unit.

We now claim that if $p_1,p_2\in M_k(B)$ are two projections
such that $\sup\{|\tau'(p_1-p_2)|: \tau'\in T(B)\}<d,$
then, for sufficiently large $m$,
$$
\sup\{|\tau(\psi_m(p_1-p_2))|:
\tau\in T(A)\}<d.
$$
To see this, assuming that $\tau_m\in T(A)$ satisfies
$|\tau_m(\psi_m(p_1-p_2))|\geq d$,
let $s$ be a weak limit of $\{\tau_m\circ \psi_m\}.$
As above, $s$ is a tracial state on $B.$ It follows that
$$
|s(p_1-p_2)|\geq d,
$$
which is a contradiction.

Fix a projection $p\in M_k(A)$ and $\ep>0.$
There is a natural number $N_1$ such that
\[ \limsup_{m\rightarrow\infty}
|\tau(\psi_m\circ \phi_n(p))-\tau(p)|<\ep/2 \]
holds for all $\tau\in T(A)$ and for all $n\geq N_1$.

Let $q\in M_{\infty}(B)$ so that $\tau'(q)=H(\hat{p})(\tau')$
for all $\tau'\in T(B).$
There is an integer $N_2>0$ such that
$$
|\tau'(q)-\tau'(\phi_n(p))|<\ep/8\rforal \tau'\in T(B)
$$
and $n\ge N_2.$
We may also assume that
there are projections $p_n\in M_k(B)$ such that
$$
\|\phi_n(p)-p_n\|<\ep/8\rforal n\ge N_2.
$$
In particular,
$$
|\tau'(\phi_n(p))-\tau'(p_n)|<\ep/8\rforal \tau'\in T(B).
$$
Hence, if $n\ge N_2,$
\[ |\tau'(q)-\tau'(p_n)|<\ep/4\rforal \tau'\in T(B). \]
Put $n=\max\{N_1,N_2\}$. It follows from the claim that there
exists a natural number $N_3$ greater than $N_1$ such that
\[ |\tau(\psi_m(q)-\psi_m(p_n))|<\ep/4\rforal \tau\in T(A)\]
and $m\geq N_3$.

Then, for sufficiently large $m$, we have
\begin{align*}
&|\tau(\psi_m(q))-\tau(p)| \\
&=|\tau(\psi_m(q))-\tau(\psi_m(p_n))|+
|\tau(\psi_m(p_n))-\tau(\psi_m(\phi_n(p)))|+
|\tau(\psi_m(\phi_n(p)))-\tau(p)| \\
&<\ep/4+\ep/8+\ep/2<\ep
\end{align*}
for all $\tau\in T(A)$, which means that
\[ \hat{p}=H'(\hat{q})=H'H(\hat{p}).\]
It follows that $H'\circ H={\rm id}_{\rho_A(K_0(A))}.$
Similarly $H\circ H'={\rm id}_{\rho_B(K_0(B))}.$
\end{proof}



\begin{Prop}\label{IIIP3}
If $(j_{\alpha*})(K_0(C(X))_+)=K_0(A_\alpha)_+$ and
$(j_{\beta*})(K_0(C(X))_+)=K_0(A_\beta)_+$,
then $\alpha$ and $\beta$ are approximately
$\tau$-conjugate if and only if they are weakly approximately
conjugate and in addition, for any  projection $f\in M_k(C(X))$
there exists $N$ such that
$$
\tau'(f\circ\gamma_n)=\tau'(f\circ\gamma_{n+1}), \ \
\tau(f\circ\sigma_n)=\tau(f\circ\sigma_{n+1}),
$$
$$
\tau(f\circ\gamma_n\circ\sigma_m)=\tau(f)
\andeqn \tau'(f\circ\sigma_n\circ\gamma_m)=\tau'(f)
$$
for every $\alpha$-invariant measure $\tau,$ $\beta$-invariant
measure $\tau'$ and all $m\ge n\ge N.$
\end{Prop}

\begin{proof}
The `only if' part is obvious. Let us prove the `if' part.


Let $f$ be an $\alpha$-projection valued function.
By Lemma \ref{IIIPpv}, there is a projection $p\in M_\infty(C(X))$
such that $\tau(f)=\tau(p)$ for every $\tau\in T(A_\alpha)$.
By the assumption, there exists $N$ such that
$$
\tau'(p\circ\gamma_n)=\tau'(p\circ\gamma_{n+1})
$$
and
$$
\tau(p\circ\gamma_n\circ\sigma_m)=\tau(p)
$$
for every $\alpha$-invariant measure $\tau,$ $\beta$-invariant
measure $\tau'$ and all $m\ge n\ge N.$
Put $q=p\circ\gamma_n$. Then $q$ is a projection of $M_\infty(C(X))$,
and so is $\beta$-projection valued. Since
\begin{align*}
& |\tau'(f\circ\gamma_n)-\tau'(q)| \\
&\leq|\tau'(f\circ\gamma_n-\phi_n(f))|+
|\tau'(\phi_n(f))-\tau'(\phi_n(p))|+|\tau'(\phi_n(p)-p\circ\gamma_n)|
+|\tau'(p\circ\gamma_n)-\tau'(q)|,
\end{align*}
we have
\[ \lim_{n\rightarrow\infty}|\tau'(f\circ\gamma_n)-\tau'(q)|=0 \]
uniformly on $T(A_\beta).$
Similarly we estimate that
\begin{align*}
& \sup_{\tau\in T(A_\alpha)}\limsup_{m\to\infty}
|\tau(f\circ\gamma_n\circ\sigma_m)-\tau(f)| \\
&\leq\sup_{\tau\in T(A_\alpha)}\limsup_{m\to\infty}(
|\tau(\psi_m(f\circ\gamma_n))-\tau(\psi_m(p\circ\gamma_n))| \\
&\, \ +|\tau(\psi_m(p\circ\gamma_n)-p\circ\gamma_n\circ\sigma_m)|
+|\tau(p\circ\gamma_n\circ\sigma_m)-\tau(p)|)
\end{align*}
converges to zero as $n\to\infty$.
\end{proof}

\section{Cantor minimal systems}

\begin{Lem}\label{Lplift}
Let $X$ be the Cantor set and $\alpha,\beta: X\to X$ be
minimal homeomorphisms.
Let $G\subset C(X, {\mathbb Z})$ be a finitely generated subgroup
containing $1_{C(X)}.$
Then there exists a finitely generated subgroup $F\supset G$ satisfying
the following:
If
$h: (j_{\alpha})_*(F)\to K_0(B)$ is a \hm\, for which
$$
h([1_A])=[1_B],\,h((j_{\alpha})_*(F_+))\subset K_0(B)_+
\andeqn
h(x)\not=0
$$
for any $x\in (j_{\alpha})_*(F_+)\setminus \{0\},$
then there is a homeomorphism
$\sigma: X\to X$ such that
$$
h(j_{\alpha})_*|_G=(j_{\beta}\circ \sigma^*)_*|_G.
$$
If one does not assume that,
for each $x\in (j_{\alpha})_*(F_+)\setminus \{0\},$ $h(x)\not=0,$
one can find a \hm\, $\phi: C(X)\to C(X)$ such that
$$
\phi_*|_G=h.
$$
\end{Lem}
\begin{proof}
Let $\{p_1,p_2,...,p_m\}$ be nonzero mutually orthogonal projections
in $C(X)$ such that $\sum_{i=1}^mp_i=1_{C(X)}$ and
the subgroup $F$ of $C(X,{\mathbb Z})$ generated by $\{p_1,p_2,...,p_m\}$
contains $G.$
Suppose that $h:(j_{\alpha})_*(F)\to K_0(B)$ is a \hm\, for which
$h([1_A])=[1_B]$, $h((j_{\alpha})_*(F_+))\subset K_0(B)_+$
and $h(x)\not=0$ if $x\in (j_{\alpha})_*(F_+)\setminus \{0\}.$
Since $\sum_{i=1}^mh([j_{\alpha}(p_i)])=[1_B]$,
by Lemma \ref{Lext3}, there are nonzero mutually orthogonal projections
$\{d_1,d_2,...,d_m\}$ in $C(X)$ such that $\sum_{i=1}^m d_i=1_{C(X)}$ and
$(j_{\beta})_*(d_i)=h([j_{\alpha}(p_i)]),$ $i=1,2,...,m.$
Since non-empty clopen subsets are homeomorphic,
there is a homeomorphism $\sigma: X\to X$ such that $\sigma^*(p_i)=d_i.$
It is then easy to verify that
$\sigma$ meets the requirements.

If we do not assume that $h(x)\not=0$ for all $x\in
(j_{\alpha})_*(F_+)\setminus \{0\},$ then in the above, we may
assume $q_i=0$ for $0<k<i\le m.$ Therefore, we may assume that
$d_i=0$ for $0<k<i\le m.$ We then fix a point $\xi_i$ in the
clopen set corresponding to $p_i$ ($i\le k$) and define
$\phi(f)=\sum_{i=1}^k f(\xi_i)d_i$ for all $f\in C(X).$ One then
sees that $\phi$ meets the requirements.
\end{proof}

\begin{Lem}\label{ILe}
Let $A=C(X)\times_{\alpha} {\mathbb Z}$ and $B=C(X)\times_{\beta}{\mathbb Z},$
where $\alpha$  and $\beta$ are minimal homeomorphisms on the Cantor set $X$.
Let $G_n$ be an increasing sequence of
finitely generated subgroups of $C(X,\Z)$ which contains $1_X$
so that $\bigcup_{n=1}(j_\alpha)_*(G_n)=K_0(A).$
Let $\lambda: K_1(A)\to K_1(B)$ be an
isomorphism. Suppose that there is a sequence of \hm s
$h_n: (j_\alpha)_*(G_n)\to K_0(B)$
for which $h_n((j_\alpha)_*((G_n)_+))\subset K_0(B)_+$
such that $h_n([1_A])=[1_B].$ Then there
exists a sequential asymptotic morphism $\{\phi_n\}: A\to B$
such that, for each $n,$ $[\phi_n]|_{(j_\alpha)_*(G_n)} $ and
$[\phi_n]|_{K_1(A)}$ are well defined
and
$$
[\phi_n]|_{(j_\alpha)_*(G_n)}=h_n \andeqn
[\phi_n]|_{K_1(A)}=\lambda
$$
Moreover, there are \hm s $\phi_n': j_{\alpha}(C(X))\to B$
such that
$$
\lim_{n\to\infty}\|\phi_n'\circ j_{\alpha}(f)-\phi_n\circ j_{\alpha}(f)\|=0
$$
for all $f\in C(X).$
\end{Lem}
\begin{proof}
We write $A=\lim_{n\to\infty}(A_n,\phi_{n,n+1})$ where $A_n$ is a
circle algebra. Let $\phi_{n,\infty}: A_n\to A$ be the \hm\,
induced by the inductive limit. It follows from 1.1 of  \cite{NT}
that we may assume that $\phi_{n, \infty}$ is injective. We may assume
that $(j_\alpha)_*(G_n)\supset (\phi_{n,\infty})_{*0}(K_0(A_n)).$
Note that
$K_1(A)={\mathbb Z}=K_1(B).$ So we may also assume that $(\phi_{n,
\infty})_{*1}(K_1(A_n))\supset K_1(A).$ Fix a finite subset ${\cal
F}\subset A.$ We may assume that ${\cal F}\subset
\phi_{n,\infty}(A_n).$ It should be noted that
$\bigcup_{n=1}^{\infty}(j_\alpha)_*((G_n)_+)=K_0(A)_+.$ Note that
$K_0(A_n)_+$ are finitely generated. We may also assume that
$$
(j_\alpha)_*((G_n)_+)\supset (\phi_{n, \infty})_*(K_0(A_n)_+).
$$
It follows from \cite{El} that there exists
a unital \hm\, $\psi_n: A_n\to B$ such that
$$
(\psi_n)_{*0}=h_n\circ(\phi_{n, \infty})_{*0}
\andeqn
(\psi_n)_{*1}=\lambda\circ(\phi_{n, \infty})_{*1},
$$
where $\lambda: K_1(A)\to K_1(B)$ is an isomorphism. Since both
$A$ and $B$ are nuclear,
by identifying
$A_n$ with $\phi_{n, \infty}(A_n),$ for any $\ep>0,$ there is a
unital contractive completely positive map $\phi_n: A\to B$ such that
$$
\|\phi_n(a)-\psi_n(a)\|<\ep\rforal a\in A_n
$$
(see 3.2 of \cite{L2}).
Note if $a, b\in A_n$ with $\|a\|,\|b\|\le 1,$
$$
\|\phi_n(ab)-\phi_n(a)\phi_n(b)\|\le \|\phi_n(ab)-\psi_n(ab)\| +
\|\psi_n(a)\psi_n(b)-\phi_n(a)\phi_n(b)\|
$$
$$
\le \ep+2\ep=3\ep.
$$
For each non-zero projection $e\in B,$ since $eBe$ is a unital a simple
non-elementary  $A\T$-algebra with real rank zero,
there is a unital  \SCA\, $C\subset eBe$ which is simple  AF-algebra
with the same ordered $K_0$-group (for example, see 2.9 of \cite{LnS}).
It is well known that there is a unital commutative \SCA\, $D\subset C$
such that $D=C(Y),$ where $Y$ is a Cantor set.
Since ${\mathbb C}^m$ as a \CA\, is semiprojective,
combinied with the above fact, one obtains a
\hm\, $\phi_n': j_{\alpha}(C(X))\to B$ such that
$$
\lim_{n\to\infty}
\|\phi_n'(j_{\alpha}(f))-\phi_n(j_{\alpha}(f))\|=0
\rforal f\in C(X).
$$
\end{proof}

\begin{Lem}\label{IL1a}
Let $A$ be a unital simple \CA\, with stable rank one and let $X$
be the Cantor set. Suppose that $h_1, h_2: C(X)\to A$ are two
monomorphisms such that $h_1(1_{C(X)})$ are both unital or both
are not unital.

Let $\ep>0$ and ${\cal F}\subset C(X)$ be a finite subset.
there is $\dt>0,$ and a finite subset ${\cal P}$ of $K_0(C(X))$
for which the following holds:
if
$$
(h_1)_*|_{\cal P}=(h_2)_*|_{\cal P},
$$
then there is a unitary $u\in A$ such that
$$
\|u^*h_1(f)u-h_2(f)\|<\ep\rforal f\in {\cal F}.
$$
\end{Lem}
\begin{proof}
Without loss of generality, we may assume that
$e=h_1(1_{C(X)})=h_2(1_{C(X)}).$ By replacing $A$ by $eAe,$ we may
further assume that both $h_1$ and $h_2$ are unital.

There exists a finite subset of mutually orthogonal projections
$\{p_1,p_2,...,p_m\}$ in $C(X)$ with $\sum_{i=1}^mp_i=1_{C(X)}$
and finitely many points $x_1,x_2,...,x_m\in X$ such that
$$
\|f-\sum_{i=1}^m f(x_i)p_i\|<\ep/2\rforal f\in {\cal F}.
$$

Let ${\cal P}=\{p_1,p_2,...,p_m\}.$ Assume
$[h_1(p_i)]=[h_2(p_i)]$ in $K_0(A).$  Since $A$ has stable rank one,
it follows that there are
partial isometries $u_i\in A$
such that
$$
u_i^*h_1(p_i)u_i=h_2(p_i), i=1,2,...,m.
$$
Since $1=\sum_{i=1}^mh_1(p_i)=\sum_{i=1}^mh_2(p_i),$
there exists a unitary $u\in A$ such that
$$
u^*h_1(p_i)u=h_2(p_i), i=1,2...,m.
$$
Then
we have
$$
\|u^*h_1(f)u-h_2(f)\|\le
\|u^*h_1(f)u-\sum_{i=1}^m f(x_i)u^*(h_1(p_i))u\|+
$$
$$
\|\sum_{i=1}^m f(x_i)u^*(h_1(p_i))u-\sum_{i=1}^m f(x_i)h_2(p_i))\|
+\|\sum_{i=1}^m f(x_i)h_2(p_i))-h_2(f)\|<\ep/2+\ep/2=\ep
$$
for all $f\in {\cal F}.$
\end{proof}

\begin{Lem}\label{IL1}
Let $A$ be a unital simple \CA\, with stable rank one  and let $X$
be the Cantor set. Suppose that $h_1, h_2: C(X)\to A$ are two
monomorphisms. Then $h_1$ and $h_2$ are approximately unitarily
equivalent if and only if $(h_1)_*=(h_2)_*$ on $K_0(C(X)).$
\end{Lem}
\begin{proof}
This is an immediate consequence of the previous lemma.
\end{proof}

\begin{Lem}\label{IL2}
Let $X$ be the Cantor set, $\alpha, \beta: X\to X$ be two minimal
homeomorphisms. Let $A=C(X)\times_{\alpha} {\mathbb Z}$ and
$B=C(X)\times_{\beta} {\mathbb Z}.$ Let $\{F_n\}$ and $\{G_n\}$ be
increasing sequences of finitely generated subgroups of $C(X,\Z)$
such that $\bigcup_{n=1}^{\infty} F_n=C(X,\Z)$ and
$\bigcup_{n=1}^{\infty}G_n =C(X,\Z),$ respectively. Suppose that
there exist sequences of
\hm s $\xi_n: (j_{\alpha})_*(F_n) \to K_0(B)$ and $\zeta_n:
(j_{\beta})_*(G_n)\to K_0(A)$ for which
$\xi_n(j_{\alpha})_*((F_n)_+)\subset K_0(B)_+$ and
$\zeta_n(j_{\beta})_*((G_n)_+)\subset K_0(A)_+$ such that
$\xi_n([1_A])=[1_B]$ and $\zeta_n([1_B])=[1_A].$
Then there exist sequential asymptotic
morphisms $\{\phi_n\}: A\to B$ and $\{\psi_n\}: B\to A$ and \hm s
$\phi_n': j_{\alpha}(C(X))\to B,$ $\psi_n': j_{\beta}(C(X))\to A,$
\hm s $\Phi_n, \Psi_n: C(X)\to C(X)$ with
$$
\lim_{n\to\infty}\|\phi_n(j_{\alpha}(f))-
\phi_n'(j_{\alpha}(f))\|=0
\andeqn
\lim_{n\to\infty}\|\phi_n(j_{\beta}(f))-
\psi_n'(j_{\beta}(f))\|=0
$$
for all $f\in C(X),$ and
\begin{eqnarray}\label{eL2}
\lim_{n\to\infty}\|\phi_n(j_{\alpha}(f))-j_{\beta}(\Phi_n(f))\|=0
\andeqn
\lim_{n\to\infty}\|\psi_n(j_{\beta}(f))-j_{\alpha}(\Psi_n(f))\|=0
\end{eqnarray}
for all $f\in C(X).$ Furthermore, if
$\xi_n((j_{\alpha})_*((F_n)_+)\setminus \{0\})\subset
K_0(B)_+\setminus \{0\}$ and $
\zeta((j_{\beta})_*((G_n)_+)\setminus \{0\})\subset
K_0(A)\setminus \{0\},$  then one can require that $\Phi_n$ and
$\Psi_n$ to be isomorphisms,
\end{Lem}
\begin{proof}
It follows from Lemma \ref{ILe} that only (\ref{eL2}) needs to be
proved. Note also as in \ref{ILe}, one has
$$
[\phi_n]|_{(j_{\alpha})_{*0}(F_n)}=\xi_n\andeqn
[\psi_n]|_{(j_{\beta})_{*0}(G_n)}=\zeta_n.
$$
By applying Lemma \ref{Lplift} and \ref{IL1a}, one obtains a
sequence of homomorphisms $\Phi_n: C(X)\to C(X)$ and a sequence of
unitaries $u_n\in B$ such that
$$
\lim_{n\to\infty}\|{\rm ad}\, u_n\circ \phi_n'\circ
j_{\alpha}(f)-j_{\beta}(\Phi_n(f))\|=0 \rforal f\in C(X).
$$
We may replace $\phi_n$ by ${\rm ad}\, u_n\circ \phi_n$  and
$\phi'$ by ${\rm ad}\, u_n\circ \phi_n'.$ We then apply the same
argument to $\psi_n.$ For the very last of the statement, we note
that, by Lemma \ref{Lplift}, under the assumption that
$\xi_n((j_{\alpha})_*((F_n)_+)\setminus \{0\})\subset
K_0(B)_+\setminus \{0\}$ and $
\zeta((j_{\beta})_*((G_n)_+)\setminus \{0\})\subset
K_0(A)\setminus \{0\},$  one can require that $\Phi_n$ and
$\Psi_n$ to be isomorphisms.
\end{proof}

\begin{Lem}\label{combinatoric}
Let $\cal P$ and $\cal Q$ be partitions of $\{1,2,\dots ,N\}$
and $\pi:{\cal P}\rightarrow {\cal Q}$ be a bijection
with $\#U=\#\pi(U)$
for all $U\in{\cal P}$.
If, for every non-empty proper subset $\cal F$ of $\cal P$,
we have
\[ \bigcup_{U\in {\cal F}}U\neq
\bigcup_{U\in {\cal F}}\pi(U), \]
then there exists a cyclic permutation $\sigma$ on
$\{1,2,\dots ,N\}$ of order $N$ such that $\sigma(U)=\pi(U)$
for all $U\in\cal P$, that is,
$\sigma(i)\in\pi(U)$ whenever $i\in U\in\cal P$.
\end{Lem}

\begin{proof}
It is easy to find a permutation $\sigma$ satisfying
$\sigma(U)=\pi(U)$ for all $U\in\cal P$. If $\sigma$ consists of
one cycle, then we have nothing to do. Otherwise we will modify
$\sigma$ so that it becomes a cyclic permutation. Let $C$ be a
cycle contained in $\sigma$. Put ${\cal F}=\{U\in{\cal P} \ : \
U\cap C\neq\emptyset\}$. If $C=\bigcup_{U\in{\cal F}}U$, then
$\cal F$ is a proper subset of $\cal P$ and
\[ \bigcup_{U\in{\cal F}}U=C=\sigma(C)=
\bigcup_{U\in{\cal F}}\pi(U), \]
which is a contradiction.
Therefore we can find $U\in\cal P$ such that $U\cap C\neq \emptyset$
and $U$ is not contained in $C$.
Take $i\in U\cap C$ and $j\in U\setminus C$.
By sending $i$ to $\sigma(j)$ and $j$ to $\sigma(i)$,
we obtain a new $\sigma$.
Then the number of its cycles is less than that of the original one.
By repeating this, we get finally a cyclic permutation $\sigma$.
\end{proof}

\begin{Lem}\label{keySapp}
Let $X$ be the Cantor set and $\alpha$, $\beta$ be minimal
homeomorphisms. Let $\cal P$ be a clopen partition of $X$. If
$[1_U]=[1_{\beta(U)}]$ in $K^0(X,\alpha)$ for all $U\in \cal P$,
then we can find a homeomorphism $\sigma\in[[\alpha]]$ such that
$\sigma\alpha\sigma^{-1}(U)=\beta(U)$ for all $U\in \cal P$.
\end{Lem}
\begin{proof}
Let $\cal Q$ denote the clopen partition
$\{\beta(U):U\in\cal P\}$.
We use the Bratteli-Vershik model for $\alpha$
(see \cite[Theorem 4.2]{HPS}).
Let
\[ {\cal P}_n=\{E(n,v,k):v\in V_n, \
k=1,2,\dots ,h(v)\} \] be a sequence of Kakutani-Rohlin partitions
such that ${\cal P}_{n+1}$ is finer than ${\cal P}_n$ for all
$n\in \N$, $\bigcup_n{\cal P}_n$ generate the topology of $X$ and
the roof sets
\[ U_n=\bigcup_{v\in V_n}E(n,v,h(v))\]
shrink to a single point in $X$. By taking a sufficiently large
$n$, we may assume that ${\cal P}_n$ is finer than $\cal P$ and
$\cal Q.$  Moreover, for every $U\in\cal P$ and $v\in V_n$,
\[ \#\{k:E(n,v,k)\subset U\}=
\#\{k:E(n,v,k)\subset \beta(U)\}\geq 1 \] since
$[1_U]=[1_{\beta(U)}]$ by the assumption. Note that these
properties hold for every $m$ larger than $n$. Because $\alpha$ is
minimal, we can find $m>n$ so that for every $v\in V_n$,
$k\in\{1,2,\dots ,h(v)\}$ and $w\in V_m$, there exists
$j\in\{1,2,\dots ,h(w)\}$ such that $E(m,w,j)\subset E(n,v,k)$.

Define a map $p_n:{\cal P}_n\rightarrow\cal P$ so that
$p(U)\supset U$ for all $U\in\cal P$.
Define $q_n:{\cal P}_n\rightarrow\cal Q$ in the same way.
Let $\cal F$ be a non-trivial proper subset of $\cal P$ and
put ${\cal G}=\{\beta(U)\in{\cal Q}:U\in\cal F\}$.
Since $\beta$ has no non-trivial invariant clopen sets,
we have
\[ \bigcup_{U\in\cal F}U\neq\bigcup_{V\in\cal G}V. \]
So there exist $v\in V_n$ and an integer $k$ such that either
$p_n(E(n,v,k))$ belongs to $\cal F$ but $q_n(E(n,v,k))$ does not
belong to $\cal G$, or $p_n(E(n,v,k))$ does not belong to $\cal F$
but $q_n(E(n,v,k))$ belongs to $\cal G$. Hence, by the choice of
$m$, for every $w\in V_m$, we get
\[ \bigcup_{p_m(E(m,w,j))\in\cal F}E(m,w,j)\neq
\bigcup_{q_m(E(m,w,j))\in\cal G}E(m,w,j). \]
By applying Lemma \ref{combinatoric},
we obtain a cyclic permutation $\sigma_w$ of order $h(w)$
such that $q_m(E(m,w,\sigma_w(j)))=\beta(p_m(E(m,w,j)))$ holds
for all $j=1,2,\dots,h(w)$.

Choose $U_0\in\cal P$ and
$j_w\in\{1,2,\dots ,h(w)\}$ for all $w\in V_m$ such that
$p_m(E(m,w,j_w))=U_0$.
For $w\in V_m$ and $j\in\{1,2,\dots ,h(w)\}$,
put $r(w,j)=\sigma_w^j(j_w)-j$.
Define a homeomorphism $\sigma$ by
$\sigma(x)=\alpha^{r(w,k)}(x)$ for $x\in E(m,w,j)$.
It is clear that $\sigma$ belongs to $[[\alpha]]$.
One can check $\sigma\alpha\sigma^{-1}(E(m,w,k))=
E(m,w,\sigma_w(k))$ for $k\neq j_w$.
Furthermore we have
\[ \bigcup_{w\in V_m}\sigma\alpha\sigma^{-1}(E(m,w,j_w))
=\bigcup_{w\in V_m}E(m,w,\sigma_w(j_w)). \]
Consequently we get $\sigma\alpha\sigma^{-1}(U)=\beta(U)$
for all $U\in\cal P$.
\end{proof}

\begin{thm}\label{TM}
Let $X$ be the Cantor set and let $\alpha, \beta: X\to X$ be two
minimal homeomorphisms. Let $A=A_{\alpha}$ and $B=A_{\beta}.$ Then
the following are equivalent:
\begin{enumerate}
\item $\alpha$ and $\beta$ are weakly approximately  conjugate;
\item There exist two unital sequential asymptotic morphisms
$\{\phi_n\}: A\to B$ and $\{\psi_n\}: B\to A;$
\item There are two increasing sequences of finitely generated
subgroups $\{F_n\}$ and $\{G_n\}$ of $K_0(C(X))$  with
$\bigcup_{n=1}F_n=\bigcup_{n=1}^{\infty} G_n=K_0(C(X))$ and two
sequences of  \hm s $\xi_n: (j_{\alpha})_*(F_n)\to K_0(B)$ with
$\xi_n((j_{\alpha})_*((F_n)_+)\subset K_0(B)_+$ and $\zeta_n:
(j_\beta)_*(G_n)\to K_0(A)$ with $\zeta_n((j_{\beta})_*((G_n)_+)
\subset K_0(A)_+$ for which $\xi_n([1_A])=[1_B],$
$\zeta_n([1_B])=[1_A],$ and for any $x\in K_0(A)_+\setminus \{0\}$
and $y\in K_0(B)_+\setminus \{0\},$ there is $N>0$ such that
$\xi_n(x)\not=0$ and $\zeta_n(y)\not=0$ for all $n\ge N;$
\item there are two increasing sequences of finitely generated
subgroups $\{F_n\}$ and $\{G_n\}$ of $K_0(C(X))$  with
$\bigcup_{n=1}F_n=\bigcup_{n=1}^{\infty} G_n=K_0(C(X))$
and two sequences of positive
\hm s $\xi_n: (j_{\alpha})_*(F_n)\to K_0(B)$ with
$\xi_n((j_{\alpha})_*((F_n)_+)\subset K_0(B)_+$ and $\zeta_n:
(j_\beta)_*(G_n)\to K_0(A)$ with $\zeta_n((j_{\beta})_*((G_n)_+)
\subset K_0(A)_+$ for which $\xi_n([1_A])=[1_B]$ and
$\zeta_n([1_B])=[1_A],$ respectively, and two sequences of
positive isomorphisms $(\sigma_n)_*, (\gamma_n)_*: K_0(C(X))\to
K_0(C(X))$ with $(\sigma_n)_*([1_{C(X)}])=[1_{C(X)}],$
$(\gamma_n)_*([1_{C(X)}])=[1_{C(X)}]$ such that, for each $x\in
K_0(C(X)),$
$$
\xi_n\circ (j_{\alpha})_*(x)=(j_{\beta})_*\circ
(\sigma_n)_*(x)\andeqn \zeta_n\circ
(j_{\beta})_*(x)=(j_{\alpha})_*\circ (\gamma_n)_*(x)
$$
for all large $n;$
\item there exist two sequential asymptotic morphisms
$\{\phi_n\}: A\to B$ and $\{\psi_n\}: B\to A$ and \hm s
$\phi_n':j_{\alpha}(C(X))\to B$ and $ \psi_n': j_{\beta}(C(X))\to
A$ such that
$$
\lim_{n\to\infty}||\phi_n\circ j_{\alpha}(f)-\phi_n'\circ
j_{\alpha}(f)||=0 \andeqn \lim_{n\to\infty}||\psi_n\circ
j_{\beta}(f)-\psi_n'\circ j_{\beta}(f)||=0
$$
for all $f\in C(X),$ and there are unital injective \hm s
$\Phi_n, \Psi_n:C(X)\to C(X)$ such that
$$
\lim_{n\to\infty}\|\phi_n'\circ j_{\alpha}(f)-j_{\beta}\circ
\Phi_n(f)\|=0 \andeqn \lim_{n\to\infty}\|\psi_n'\circ
j_{\beta}(f)-j_{\alpha}\circ \Psi_n(f)\|=0
$$
for all $f\in C(X);$
\end{enumerate}
\end{thm}
\begin{proof}
(1) $\Rightarrow$ (2) follows from Proposition \ref{wac>asym}.

(2) $\Rightarrow$ (3) follows from Definition \ref{DKK}. But note
that, when $A$ (and $B$) are simple, by Lemma \ref{ILeffchoi} and by
Definition \ref{DKK}, for any nonzero projection $p,$
$[\phi_n]([p])$ can not be zero in the simple dimension group
$K_0(B)_+.$

(3) $\Rightarrow$ (4) follows from Lemma \ref{Lplift}.

(3) $\Rightarrow$ (5) follows from Lemma \ref{IL2}.

(4) $\Rightarrow$ (3). Note that the requirements that
$\xi_n(x)>0$ eventually follows from the fact that $K_0(A)$ is a
simple dimension group and $mx\ge [1_A]$ if $x\in
K_0(A)_+\setminus \{0\}.$

(5)$\Rightarrow$(1):
Fix $\ep>0$ and a finite subset ${\cal F}\subset C(X).$
Let ${\cal P}$ be a clopen partition  for $X$ and
$p_1,p_2,...,p_m$ be corresponding projections in $C(X).$
We choose ${\cal P}$ so that
$$
\|f-\sum_{i=1}^m f(x_i)p_i\|<\ep/2\rforal f\in {\cal F},
$$
where $x_i$ in the clopen set in ${\cal P}$ associated with $p_i.$
Let ${\cal G}={\cal F}\cup\{p_1,...,p_m\}.$
Without loss of generality, we may assume
that
\begin{eqnarray}\label{eP}
\|\phi_n'\circ j_{\alpha}(f)-j_{\beta}\circ \Phi_n(f)\|<\ep/4
\end{eqnarray}
for all $f\in {\cal G}\cup \alpha^*({\cal G}).$
We have
\begin{eqnarray}\label{eP2}
\lim_{n\to\infty}\|\phi_n(u_{\alpha})\phi_n'(j_{\alpha}(f))
\phi_n(u_{\alpha}^*)
-\phi_n'(j_{\alpha}(\alpha^*(f)))\|=0
\end{eqnarray}
for all $f\in C(X).$
There are unitaries $z_n\in B$ such that
\begin{eqnarray}\label{eP3}
\lim_{n\to\infty}\|\phi_n(u_{\alpha})-z_n\|=0.
\end{eqnarray}
It follows from this and (\ref{eP}) and (\ref{eP2}) that
\begin{eqnarray}
\|z_nj_{\beta}(\Phi_n(p_i))z_n^*- j_{\beta}
(\Phi_n(\alpha^*(p_i)))\|<\ep/4, \, i=1,2,...,m.
\end{eqnarray}
In particular, (for $0<\ep<1/2$)
$$
[j_{\beta}\circ \Phi_n(p_i)]=[j_{\beta}\circ
\Phi_n(\alpha^*(p_i))],i=1,2,...,m.
$$
This implies, by Lemma \ref{keySapp},
there is a homeomorphism $\sigma\in [[\beta]]$ such
that
$$
(\sigma\beta\sigma^{-1})^*(\Phi_n(p_i))=\Phi_n(\alpha^*(p_i)),
 \ \ i=1,2,...,m.
$$
It follows that
$$
\|(\sigma\beta\sigma^{-1})\circ \Phi_n(f)-\Phi_n\circ
\alpha^*(f)\|<\ep\rforal f\in {\cal F}.
$$
Thus we obtain a homeomorphism $\sigma': X\to X$
such that
$$
\|(\sigma'\beta(\sigma')^{-1})^*(f)-\alpha^*(f)\|<
\ep\rforal f\in {\cal F}.
$$
The same argument shows that
there exists a homeomorphism $\gamma: X\to X$ such that
$$
\|(\gamma\alpha\gamma^{-1})^*(f)-\beta^*(f)\|<\ep \rforal f\in {\cal F}.
$$
\end{proof}

\begin{Remark}\label{RTM}
{\rm In the proof of (5) $\Rightarrow$ (1),  we remark that for a
given finitely generated subgroup of $K_0(A) $ (and $K_0(B)$),
$\sigma'$ (and $\gamma$) may be chosen so it induces the same map
on that subgroup as that of $\Phi_n$ (and $\Psi_n$), since
$\sigma\in [[\alpha]].$}
\end{Remark}

\begin{Cor}\label{IIC1}
Let $X$ be the Cantor set and $\alpha, \beta: X\to X$
 be two minimal homeomorphisms. Suppose that
 there is a positive \hm\, $h: K^0\xa\to K^0\xb$
 with $h([1_A])=[1_B].$ Then
 there exists a sequence of homeomorphisms $\sigma_n: X\to X$
 such that
 $$
  \lim_{n\to\infty}(\sigma_n\beta\sigma_n^{-1})^*(f)=\alpha^*(f)\rforal
  f\in C(X).
  $$
  Moreover, there exists a sequence of $\gamma_n\in [[\beta]]$
 and a unital injective \hm\, $\Phi: C(X)\to C(X)$ such that
$$
\lim_{n\to\infty}(\gamma_n\beta\gamma_n^{-1})^*\circ \Phi(f)=
\Phi\circ \alpha^*(f)\rforal f \in C(X).
$$
Moreover, $(j_{\beta})_{*0}\circ \Phi_{*0}=h\circ (j_{\alpha})_{*0}.$
\end{Cor}

  \begin{proof}
Since both $K^0\xa$ and $K^0\xb$ are simple dimension groups,
we conclude that
$$
h\circ (j_{\alpha})_{*0}(C(X,\Z)_+\setminus \{0\})\subset
K^0\xb_+\setminus \{0\}.
$$
It follows from Theorem \ref{PreT} that
there is an injective  unital positive
\hm\, $\dt: C(X,\Z)\to C(X,\Z)$ such that
$$
(j_{\beta})_{*0}\circ \dt=h\circ(j_{\alpha})_{*0}.
$$
One obtains an injective unital \hm\, $\Phi: C(X)\to C(X)$ such that
$\Phi_{*0}=\dt.$
In the proof of Theorem \ref{TM}, one can choose $\Phi_n=\Phi$.
\end{proof}

\begin{Def}\label{Ddiv}
Let $(G, G_+,u)$ be a scaled ordered abelian group. Set
$$
D(G,u)=\{n\in \N: ne=u\,\,\,{\rm for\,\,\, some\,\,\,} e\in G_+\}.
$$
\end{Def}

\begin{Def}
Let $\xa$ be a Cantor minimal system.
By the periodic spectrum of $\xa$ or $\alpha$, we mean
the set of natural numbers $p$ for which there are disjoint
clopen sets $U,\alpha(U),\dots,\alpha^{p-1}(U)$ whose union
is $X$.
\end{Def}
It is obvious that  $p$ is a periodic spectrum if and only if
$p\in D(K_0(A_{\alpha}), [1_{A_{\alpha}}]).$


\begin{thm}\label{wac}
Let $X$ be the Cantor set and $\alpha, \beta: X\to X$ be two
minimal homeomorphisms. Then the following are equivalent.
\begin{enumerate}
\item $\alpha$ and $\beta$ are weakly approximately conjugate.
\item The periodic spectrum of $\alpha$ and $\beta$ agree.
\item $D(K_0(A_{\alpha}),[1_{A_{\alpha}}])=D(K_0(A_{\beta}), [1_{A_{\beta}}]).$

\end{enumerate}

\end{thm}
\begin{proof}
We have seen (2) $\Leftrightarrow$ (3).

(1) $\Rightarrow$ (3): Suppose that $p\in D(K_0(A_{\alpha}),
[1_{A_{\alpha}}]).$ Then there is $e\in K_0(A_{\alpha})_+$ such
that $pe=[1_{A_{\alpha}}].$ Let $G\subset K_0(A_{\alpha})$
generated by $e.$ Suppose that $\{\psi_n\}: A_{\alpha}\to
A_{\beta}$ is a unital asymptotic morphism given by Proposition
\ref{wac>asym}. It follows that there are \hm s $h_n: G\to
K_0(A_{\beta})$ such that $h_n(e)\ge 0,$  $h_n(e)\not=0$ and
$h_n([1_{A_{\alpha}}])=[1_{A_{\beta}}].$
Thus $ph_n(e)=[1_{\beta}].$ This
proves that $D(K_0(A_{\alpha}),[1_{A_{\alpha}}])\subset
D(K_0(A_{\beta}), [1_{A_{\beta}}]).$ Applying the same argument to
$A_{\beta},$ one obtains
$D(K_0(A_{\alpha}),[1_{A_{\alpha}}])=D(K_0(A_{\beta}),
[1_{A_{\beta}}]).$

(3) $\Rightarrow$ (1): Write $A=\lim_{n\to\infty}(A_n, \phi_n),$
where $A_n$ is a circle algebra and $\phi_n$ are injective. Thus
$K_0(A)=\lim_{k\to\infty} (K_0(A_n),(\phi_n)_*)$ where each
$K_0(A_n)=\Z^{r_n}$ with the usual order. Fix $n,$ we may assume
that there is $u_n\in K_0(A_n)_+$ such that
$\phi_{n,\infty}(u_n)=[1_A].$ Denote $e_1=(1,0,...,0),$
$e_2=(0,1,0,...,0),...,e_{r_n}=(0,...,0,1).$ We may assume that
$u_n=\sum_{i=1}^{r_n}m_ie_i,$ where $m_i>0.$ Let $p$ be the
greatest common divisor of $m_1,m_2,...,m_{r_n}.$ Write
$m_i=pk_i,$ where $k_i\in \N$ and $k_1,k_2,...,k_{r_n}$
are relatively prime. Choose an integer $N>0$ such that
$$
\{n\in \N : n\ge N\}\subset \N k_1+\N k_2+\cdots +\N k_{r_n}.
$$

 Write $B=\lim_{n\to\infty}(B_n,\psi_n),$
 where each $B_n$ is a circle algebra and $\psi_n$ are
 injective. Thus $K_0(B)=\lim_{n\to\infty} (K_0(B_n), (\psi_n)_*),$
where $K_0(B_n)=\Z^{s_n}$ with the usual order.  We may assume
that $v_n\in K_0(B_n)_+$ such that $\psi_{n,\infty}(v_n)=[1_B].$
Since $K_0(B)$ is a simple dimension group and $p$ divides $[1_B]$
in $K_0(B),$ we may assume that
$$
[1_B]=(pd_1,pd_2,...,pd_{s_n})\,\,\,{\rm with}\,\,\,
d_i>N,\,i=1,2,...,s_n.
$$
Thus there are $N(i,j)\in \N$ such that
$$
d_i=\sum_{j=1}^{r_n} N(i,j)k_i\,\,\,i=1,2,...,s_n
$$
(in $K_0(B_n)$).  Define a $s_n\times r_n$ matrix
$T_n=(N(i,j))_{s_n\times r_n}.$ Then $T_n$ gives a positive \hm\,
from $\Z^{r_n}$ to $\Z^{s_n}$ which maps $u_n$ to $v_n.$ Note that
we may write $A_n=\bigoplus_{i=1}^{r_n}C(S^1, M_{m_i}).$
Using a
point-evaluation, we obtain a unital \hm\, $\theta_n: \phi_{n,
\infty}(A_n)\to \bigoplus_{i=1}^{r_n}M_{m_i}.$ One checks that
$(\theta_n)_{*0}={\rm id}_{K_0(A_n)}.$ Write $B_n=\bigoplus_{j=1}^{s_n}
C(S^1, M_{pd_j}).$
 There is a unital \hm\,
$h_n:\bigoplus_{i=1}^{r_n}M_{m_i}\to \bigoplus_{j=1}^{s_n}M_{pd_j}$ such
that $(h_n)_{*0}=T.$ It follows from 5.2 of \cite{P} (see also
3.2 (1) of \cite{L2}) that there is a
\morp\, $L_n: A\to \bigoplus_{j=1}^{s_n}M_{pd_j}$ such that
$(L_n)|_{\phi_{n,\infty}(A_n)}=h_n\circ \theta_n.$ Now define
$\Phi_n=\psi_{n,\infty}\circ L_n.$ It follows that $\{\Phi_n\}$ is
a unital sequential asymptotic morphism from $A$ to $B.$

The same argument allows one to construct a unital sequential
asymptotic morphism from $B$ to $A.$ It follows from Theorem \ref{TM} that
$\alpha$ and $\beta$ are weakly approximately conjugate.
\end{proof}

\begin{Remark}
There is a direct dynamical proof of the implication 
from (2) (or (3)) to (1) in the theorem above.
See \cite{M}.
\end{Remark}

\section{Other approximate conjugacy in the Cantor system}


\begin{Lem}\label{IVLlift}
Let $F_1, F_2$ and $G$ be countable abelian torsion free groups
with $g\in G$ and $f_i\in F_i$ ($i=1,2$). Suppose that there is a
surjective \hm\, $\pi_i: F_i\to G$ such that $\pi_i(f_i)=g$
($i=1,2$). Suppose that $D(F_1, f_1)=D(F_2, f_2).$ Then, for each
finitely generated group $G_1\subset F_1$ containing $f_1,$ there
exists a \hm\, $h: G_1\to F_2$ such that $\pi_2\circ
h=\pi_1|_{G_1}$ and $h(f_1)=f_2.$
\end{Lem}

\begin{proof}
Write $G_1=\Z\oplus \Z\oplus\cdots \oplus \Z$  (there  are $r$
copies of $\Z$). Suppose that $f_1=(m_1,m_2,\dots, m_r),$ where
$m_i\in \Z.$  Let $p$ be the greatest common divisor of
$m_1,m_2,...,m_r.$ We may write $m_i=pk_i,$ $i=1,2,...,r$ and
$k_1,k_2,...,k_r$ are relatively prime. By the assumption there is
$f_0\in F_2$ such that $pf_0=f_2.$ Denote $e_1=(1,0,...,0),$
$e_2=(0,1,0,...,0),...,e_r=(0,...,0,1)$ in $G_1.$ Since
$k_1,k_2,...,k_r$ are relatively prime, there are
$n_1,n_2,...,n_r\in \Z$ such that
$$
\sum_{i=1}^rk_in_i=1.
$$
On the other hand, there are $g_i\in F_2$ such that
$\pi_2(g_i)=\pi_1(e_i),$ $i=1,2,...,r.$ Since $G$ is torsion free,
$f_{00}=\sum_{i=1}^rk_ig_i-f_0\in {\rm ker}\pi_2.$

Put $x_i=n_if_{00},$ $i=1,2,...,r.$ Define $h: G_1\to F_2$ by
$h(e_i)=g_i-x_i,$ $i=1,2,...,r.$ Then $\pi_2\circ h=\pi_1$
and
$$
h(f_1)=p\sum_{i=1}^rk_i(g_i-x_i)=p[(\sum_{i=1}^rk_ig_i)-f_{00}]=pf_0=f_2.
$$
\end{proof}

\begin{thm}\label{IIITwt=orb}
Let $X$ be the Cantor set and $\alpha, \beta: X\to X$ be two
minimal homeomorphisms. The following are equivalent.
\begin{enumerate}
\item $\alpha$ and $\beta$ are approximately
$\tau$-conjugate.
\item $\alpha$ and $\beta$ are orbit equivalent and have
the same periodic spectrum.
\item $\rho_A(K_0(A_\alpha))$ is unital order isomorphic to 
$\rho_B(K_0(A_\beta))$
and $D(K_0(A_{\alpha}), [1_{A_{\alpha}}])=D(K_0(A_{\beta}),
[1_{A_{\beta}}]).$
\item There exist unital sequential asymptotic morphisms
$\phi_n:A_{\alpha}\rightarrow A_{\beta}$ and
$\psi_n:A_{\beta}\rightarrow A_{\alpha}$ such that the following
are satisfied: for any projections $p\in A$ and $q\in B$ there
exists $N$ such that
\begin{eqnarray}\label{etr=10}
\tau([\psi_n]([q]))=\tau([\psi_{n+1}]([q])), \,\,\,
\tau'([\phi_n]([p]))=\tau'([\phi_{n+1}]([p])),
\end{eqnarray}
\begin{eqnarray}\label{ett=1}
\tau([\psi_m\circ \phi_n]([p]))=\tau([p]) \andeqn
\tau'([\phi_m\circ\psi_n]([q]))=\tau'([q])
\end{eqnarray}
 for every trace
$\tau$ on $A,$ every trace $\tau'$ on $B$ and all $m\ge n\ge N.$
\end{enumerate}
\end{thm}
\begin{proof}
(1) $\Rightarrow$ (4): This follows immediately from the
definition and Proposition \ref{IIIP2}.

(4) $\Rightarrow$ (3): The fact that $D(K_0(A_{\alpha}),
[1_{A_{\alpha}}])=D(K_0(A_{\beta}),[1_{A_{\beta}}])$ follows from
Theorem \ref{wac}. Then one applies Theorem \ref{IIITtr}.

(3) $\Rightarrow$ (1): Let $A=A_{\alpha}$ and $B=A_{\beta}$ and
let $\{G_n\}$ be an increasing sequence of finitely generated
subgroups of $C(X, \Z)$ containing $[1_{C(X)}]$ such that
$\bigcup_{n=1}^{\infty} G_n=C(X,\Z).$ Let $s:\rho_A(K_0(A))\to
\rho_B(K_0(B))$ be the order isomorphism preserving the unit. It
follows from Lemma \ref{IVLlift}, there exists a \hm\, $h_n: G_n\to
K_0(B)$ which maps $[1_{C(X)}]$ to $[1_B]$ and $\rho_B\circ
h_n=(s\circ \rho_A\circ (j_{\alpha})_{*0})|_{G_n},$ $n=1,2,....$
Then $h_n$ is a positive \hm\, which maps $[1_{C(X)}]$ to $[1_B].$
It follows from Theorem \ref{TM} ((2) $\Rightarrow$ (1)) that
there exists a sequence of homeomorphisms $\sigma_n: X\to X$ such
that
$$
\lim_{n\to\infty}(\sigma_n\circ \alpha\circ \sigma_n^{-1})^*(f)
=\beta^*(f)\rforal f\in C(X).
$$
Similarly, we obtain a sequence of homeomorphisms $\gamma_n: X\to
X$ such that $$ \lim_{n\to\infty} (\gamma_n\circ \beta\circ
\gamma_n^{-1})^*(f)=\alpha^*(f)\rforal f\in C(X).
$$
So $\alpha$ and $\beta$ are weakly approximately  conjugate.
 Moreover, from the
definition of $\phi_n$ and Remark \ref{RTM}  we see that they are
actually  approximately  $\tau$-conjugate. For example, if $p\in
C(X)$ is a projection, then, by the proof of Theorem \ref{TM},
$$
[j_{\beta}\circ \Phi_n(p)]=h_n\circ s\circ \rho_A\circ (j_{\alpha})_*([p])
$$
for all large $n.$
Then
$$
\rho_B([j_{\beta}\circ \Phi_n(p)])=s\circ \rho_A\circ (j_{\alpha})_*([p]),
$$
is a constant (for sufficiently large $n$).


(2) $\Leftrightarrow$ (3). This follows from Theorem 2.3 of \cite{GPS}.

\end{proof}

There is a general version of $K$-conjugacy. However, we only introduce
this notion for the Cantor set for simplicity.

\begin{Def}\label{DKC}
Let $X$ be the Cantor set, $\alpha,\beta: X\to X$ be two minimal
homeomorphisms. We say $\alpha$ and $\beta$ are {\it
approximately $K$-conjugate} if there are homeomorphisms
$\sigma_n, \gamma_n: X\to X$ such that
\[
\lim_{n\to\infty}(\sigma_n\alpha\sigma_n^{-1})^*(f)=\beta^*(f)
\andeqn
\lim_{n\to\infty}(\gamma_n\beta\gamma_n^{-1})^*(f)=\alpha^*(f)\rforal
f\in C(X),
\]
in addition, for any $g\in C(X,\Z),$ there exists $N$ such that
$$
(j_{\alpha})_{*0}(g\circ \sigma_n)=
(j_{\alpha})_{*0}(g\circ\sigma_{n+1}), \,\,\,
(j_{\beta})_{*0}(g\circ \gamma_n)=(j_{\beta})_{*0}(g\circ \gamma_{n+1}),
$$
$$
(j_{\alpha})_{*0}(g\circ \gamma_n\circ \sigma_{m})=(j_{\alpha})_{*0}(g) \andeqn
(j_{\beta})_{*0}(g\circ \sigma_n\circ \gamma_m)=(j_{\beta})_{*0}(g)
$$
for all $m\ge n\ge N.$
\end{Def}

\begin{thm}\label{Tapp=str}
Let $X$ be the Cantor set and
$\alpha, \beta: X\to X$ be minimal homeomorphisms.
Let $A=A_{\alpha}$ and $B=A_{\beta}.$
Then the following are equivalent.
\begin{enumerate}
\item $\alpha$ and $\beta$ are  approximately $K$-conjugate.
\item $C(X)\times_{\alpha}\Z\cong C(X)\times_{\beta}\Z.$
\item $(K^0(X,\alpha),K^0(X, \alpha)_+,[1])$ is unital order isomorphic to
$(K^0(X, \beta), K^0(X, \beta)_+, [1])$.
\item
There exists a homeomorphism $F: X\to X,$ and exist sequences
$\sigma_n\in [[\alpha]]$ and $\gamma_n\in [[\beta]]$ such that
$$
\lim_{n\to\infty}
(F^{-1}\circ\sigma_n\circ \alpha\circ \sigma_n^{-1}\circ F)^*(f)
=\beta^*(f)\andeqn
\lim_{n\to\infty} (F\circ \gamma_n\circ\beta\circ \gamma_n^{-1}\circ F^{-1})^*
(f)=\alpha^*(f)
$$
for all $f\in C(X).$
\item $\alpha$ and $\beta$ are strongly orbit equivalent.
\item There exist asymptotic morphisms $\phi_n:A\rightarrow B$
and $\psi_n:B\rightarrow A$ such that the following are satisfied:
there is a unit preserving  order isomorphism  $\Phi:
(K_0(A), K_0(A)_+, [1_A])\to (K_0(B), K_0(B)_+, [1_A])$
such that, for any finite subsets ${\cal G}\subset K_0(A)$
and ${\cal F}\subset K_0(B),$ there exists $N>0$ for which
$$
[\phi_n]|_{\cal G}=\Phi|_{\cal G}\andeqn [\psi_n]|_{\cal F}
=\Phi^{-1}|_{\cal F},
$$
whenever $n\ge N.$
\end{enumerate}
\end{thm}

\begin{proof}
(1) $\Rightarrow$ (6):
 This follows from the definition of approximate
$K$-conjugacy and Proposition \ref{wac>asym}.

(6) $\Rightarrow$ (3) is obvious.

(3) $\Rightarrow$ (4): This follows from Theorem \ref{PreT} and
the proof of Corollary \ref{IIC1}. Note this time the map $\Phi$
in Corollary \ref{IIC1} can be chosen to be isomorphism.

(4) $\Rightarrow$ (1): This is obvious.

(5) $\Leftrightarrow$ (2) follows from \cite{GPS}.

(2) $\Leftrightarrow$ (3) follows from \cite{El} (see also
\cite{GPS}).

\end{proof}

\begin{Remark}\label{RF}
Let us consider the case $X=S^1.$ It is well-known that every
minimal homeomorphism on $S^1$ is conjugate to an irrational
rotation $R_\alpha$ (see \cite[Proposition 6.4]{PY} and \cite[Ex. 11.2.4, p. 400]{KH},
for example). In this
case, when are two irrational rotations $R_\alpha$ and $R_\beta$
weakly approximately conjugate? By the following easy observation,
we can see that they are weakly approximately conjugate if and
only if they are actually conjugate.

We identify $S^1$ with $\T=\R/\Z$. Let us denote the set of
orientation preserving homeomorphisms on $\T$ by $\Homeo^+(\T)$.
For $\phi\in\Homeo^+(\T)$ we denote the rotation number of $\phi$
by $r(\phi)\in\T$. The reader may refer to \cite[Section 6]{PY} or
\cite[Chapter 11]{KH} for the definition and the elementary
property of $r(\phi)$.

By \cite[11.1.3]{KH}, $r(\sigma\phi\sigma^{-1})=r(\phi)$
for all $\phi,\sigma\in\Homeo^+(\T)$.
Furthermore, by \cite[11.1.6]{KH},
the map $r:\Homeo^+(\T)\rightarrow\T$ is continuous,
where $\Homeo^+(\T)$ is equipped with the so called uniform topology
(see p.710 of \cite{KH}).
It is not hard to see that the uniform topology of $\Homeo(X)$
agrees with the topology induced from pointwise convergence
in norm on $C(X)$ when $X$ is compact and metrizable.

Suppose that there exists a sequence of homeomorphisms $\sigma_n$
such that $\sigma_nR_\alpha\sigma_n^{-1}$ converges to $R_\beta$.
We may assume that $\sigma_n$ belongs to $\Homeo^+(\T)$ for all
$n\in\N$ or $\sigma_n\circ \pi$ belongs to $\Homeo^+(\T)$ for all
$n\in\N$, where $\pi$ is given by $\pi(t)=-t$. When
$\sigma_n\in\Homeo^+(\T)$, we have
\[ \beta=r(R_\beta)=
\lim_{n\rightarrow\infty}r(\sigma_nR_\alpha\sigma_n^{-1})=
\lim_{n\rightarrow\infty}r(R_\alpha)=\alpha. \]
Thus, $R_\beta$ is conjugate to $R_\alpha$.
If $\sigma_n\pi\in\Homeo^+(\T)$,
\[ \beta=r(R_\beta)=
\lim_{n\rightarrow\infty}r(\sigma_n\pi R_{-\alpha}\pi\sigma_n^{-1})=
\lim_{n\rightarrow\infty}r(R_{-\alpha})=-\alpha, \]
which also means that $R_\beta$ is conjugate to $R_\alpha$.

\vspace{0.2in}

Approximate conjugacy for $X=S^1\times Y,$ where 
$Y$ is the Cantor set, will be discussed in a subsequent paper.  
\end{Remark}

\end{document}